\documentclass[10pt,sort&compress,3p,letterpaper]{elsarticle}

\usepackage[T1]{fontenc}
\usepackage[bitstream-charter]{mathdesign}
\usepackage[scaled]{berasans}

\usepackage{geometry}
\usepackage[reqno]{amsmath} 
\usepackage{mathtools}
\usepackage{bm}

\usepackage{xspace}                   
\usepackage{xifthen}
\usepackage{ifdraft}
\usepackage{paralist}                 
\usepackage[modulo]{lineno}
\usepackage{soul}                     
\usepackage[notref,notcite,final]{showkeys}

\usepackage{siunitx}                 
\usepackage[usenames,dvipsnames,svgnames,table]{xcolor}

\expandafter\let\csname ver@natbib.sty\endcsname\relax
\expandafter\let\csname c@author\endcsname\relax
\usepackage[
  backend=bibtex,
  style=alphabetic,
  maxbibnames=99,
  maxcitenames=2,
  sorting=nyt,
  url=false,
  isbn=false,
  firstinits=true,
  sortcites
]{biblatex}

\renewcommand\citet[1]{\textcite{#1}}
\definecolor{ref-darkblue}{rgb}{0.03,0.3,0.62}
\definecolor{ref-darkorange}{rgb}{1,0.55,0}
\definecolor{ref-turquoise}{rgb}{0.25,0.88,0.82}

\usepackage[
  colorlinks,
  pdfauthor   = {Abtin Rahimian},
  pdftitle    = {Contact aware vesicle simulation},
  pdfcreator  = {pdflatex},
  pdfpagemode = UseNone
]{hyperref}

\usepackage[capitalize]{cleveref}
\crefname{equation}{Eq.}{Eqs.}
\Crefname{equation}{Equation}{Equations}
\crefname{algocf}{Alg.}{Algs.}
\Crefname{algocf}{Algorithm}{Algorithms}
\crefrangelabelformat{equation}{(#3#1#4--#5#2#6)}
\crefmultiformat{equation}{Eqs.~(#2#1#3)}{ and~(#2#1#3)}{, (#2#1#3)}{, and~(#2#1#3)}
\Crefmultiformat{equation}{Equations~(#2#1#3)}{ and~(#2#1#3)}{, (#2#1#3)}{, and~(#2#1#3)}

\usepackage{relsize}
\def\abbrevsize{.85}
\newcommand\abbrevformat[1]{\textscale{\abbrevsize}{#1}}
\newcommand\define[1]{\abbrevformat{#1}}
\newcommand\abbrev[1]{\abbrevformat{#1}}

\makeatletter
\newcommand{\hypertargetraised}[1]{\Hy@raisedlink{\hypertarget{#1}{}}}
\makeatother

\usepackage{titlesec}
\titleformat{\subsubsection}[runin]{\normalfont\normalsize\itshape}{\thesubsubsection.}{.5em}{}[.\hspace{.5em}]
\titleformat{\paragraph}[runin]{\normalfont\normalsize\itshape}{\theparagraph}{}{}[.\hspace{.5em}]
\titlespacing*{\section}{0pt}{*3}{*1.5}
\titlespacing*{\subsection}{0pt}{*2.5}{*1}

\usepackage{amsthm}

\numberwithin{equation}{section}

\newcommand\para[1]{\paragraph*{#1}}

\usepackage[linesnumbered,lined,ruled,vlined,commentsnumbered]{algorithm2e}

\usepackage[margin=10pt,font=small,textfont=it,labelfont=bf,justification=justified]{caption}
\newcommand{\algcaption}[3]{
        \ifthenelse{\isempty{#3}}
                   {\caption[#1]{{\sc #2.} \label{#1}}}
                   {\caption[#1]{{\sc #2.} \newline\small{#3} \label{#1}}}
        }

\usepackage{graphicx}
\graphicspath{{figs/}}
\DeclareGraphicsExtensions{.pdf,.png,.jpg}

\usepackage[export]{adjustbox}  
\usepackage[
  format=hang,
  singlelinecheck=false,
  font={small},
  labelfont=bf,
  justification=centering,
]{subfig}

\newcommand{\mcaption}[3]{
  \ifthenelse{\isempty{#2}}
             {\caption[#1]{#3 \label{#1}}}
             {\caption[#1]{{\sc #2.} #3 \label{#1}}}
}

\usepackage{tikz,pgfplots,import}
\usetikzlibrary{external}
\usetikzlibrary{patterns}

\newif\ifPlotTikz
\ifPlotTikz
\pgfplotsset{compat=newest}
\usepgfplotslibrary{fillbetween}
\usetikzlibrary{arrows.meta}
\usetikzlibrary{backgrounds}
\usetikzlibrary{pgfplots.groupplots}
\usetikzlibrary{plotmarks}

\pgfplotsset{plot coordinates/math parser=false}
\pgfkeys{/pgf/images/include external/.code=\includegraphics{#1}}
\newcommand{\includepgf}[2][1]{
\beginpgfgraphicnamed{#2}%
\tikzsetnextfilename{external-#2}%
\scalebox{#1}{\subimport{figs/}{#2.pgf}}%
\endpgfgraphicnamed%
}

\makeatletter
\tikzset{
    every picture/.style={
        execute at begin picture={
            \let\ref\@refstar
        }
    }
}
\makeatother
\else
\newcommand{\includepgf}[2][1]{
\scalebox{#1}{\includegraphics[]{external-#2}}%
}
\fi
\newlength\figureheight
\newlength\figurewidth

\definecolor{plt-blue}{rgb}{0.0078,0.2980,0.7961}
\definecolor{plt-orange}{rgb}{1.0000,0.6431,0.2627}
\definecolor{plt-purple}{rgb}{1.0000,0.2863,0.5255}
\definecolor{plt-violet}{rgb}{0.6118,0.1765,1.0000}

\usepackage{tabularx}     
\usepackage{booktabs}     
\usepackage{multirow}
\usepackage{collcell}     
\usepackage{dcolumn}      
\usepackage{ctable}       

\usepackage[normal]{engord} 

\let\d\undefined
\let\O\undefined
\DeclareMathOperator{\Grad}{\nabla}
\DeclareMathOperator{\Div}{\nabla\cdot}

\newcommand\defeq{\mathrel{\mathop :}=}

\newcommand\norm[2][{}]{\lVert #2 \rVert_{#1}}

\newcommand\contract[3][{\cdot}]{#2\mathop{#1}#3}
\newcommand\inner[2]{\contract{#1}{#2}}

\newcommand\d{\ensuremath{\,\mathrm{d}}}
\newcommand\dline{\ensuremath{\d s}}
\newcommand\dsurf{\ensuremath{\d A}}

\renewcommand\perp[1]{{#1}^{\bot}}

\newcommand\O[1]{\ensuremath{\mathcal{O}\left(#1\right)}}
\newcommand\ordinal[1]{\ifthenelse{\isin{#1}{abcdefghijklmnopqrstuvwxyz}}{\ensuremath{#1^\mathrm{th}}}{\engordnumber{#1}}}

\let\vector\undefined
\newcommand\mathbfsf[1]{\bm{\mathsf{#1}}}
\newcommand\discrete[1]{\mathsf{#1}}

\newcommand\scalar[1]{#1}
\newcommand\scalard[1]{\discrete{#1}}
\newcommand\vector[1]{\bm{#1}}
\newcommand\vectord[1]{\mathbfsf{#1}}
\newcommand\linop[1]{\mathbf{#1}}
\newcommand\linopd[1]{\mathbfsf{#1}}
\newcommand\conv[1]{\mathcal{#1}}

\newcommand\reynolds{\ensuremath{R\kern-.06em e}\xspace}  
\newcommand\capillary{\ensuremath{C\kern-.16em a}\xspace} 

\newcommand\vx{\vector{x}}
\newcommand\vX{\vector{X}}

\newcommand\vXd{\vectord{X}}
\newcommand\vy{\vector{y}}
\newcommand\vY{\vector{Y}}
\newcommand\vYd{\vectord{Y}}
\newcommand\vn{\vector{n}}

\newcommand\vu{\vector{u}}

\newcommand\dt{\ensuremath{\Delta t}}

\newcommand\slyr{single-layer\xspace}
\newcommand\dlyr{double-layer\xspace}
\newcommand\ns{near-singular\xspace}

\def\vel{\ensuremath{\vector{u}}}
\def\vbg{\ensuremath{\vu^\infty}}

\def\Fb{\ensuremath{\vector{f}_b}}
\def\Fs{\ensuremath{\vector{f}_\sigma}}
\def\Fc{\ensuremath{\vector{f}_c}}
\def\ten{\scalar{\sigma}}
\newcommand\vf{\vector{f}}

\newcommand\twod{{\small 2}\textsc{d}\xspace}
\newcommand\threed{{\small 3}\textsc{d}\xspace}
\def\gmres{\abbrev{GMRES}\xspace}
\def\li{\abbrev{LI}\xspace}
\def\gi{\abbrev{GI}\xspace}
\def\lic{\abbrev{CLI}\xspace}

\newcommand\sdc[1][]{\abbrev{SDC\ifthenelse{\isempty{#1}}{}{\kern 1pt}#1}\xspace}
\def\emdash/{\kern 0.2em---\kern 0.2em}

\newcommand\e[1]{\ensuremath{e{#1}}}

\newcommand\scid[2][1]{\ensuremath{\ifthenelse{\equal{#1}{1}}{}{#1\times}10^{#2}}}
\newcommand\sci[2][1]{\scid[#1]{#2}}
\newcommand\andor[3][-.15em]{#2/{\kern #1}#3}

\def\colres{collision resolution\xspace}
\newcommand\comp[1]{\bar{#1}}
\def\minsep{\ensuremath{d_{m}}\xspace}
\usetikzlibrary{fit,backgrounds}

\begin{document}

\title{Contact-aware simulations of particulate Stokesian suspensions}

\author[nyu]{Libin Lu} \ead{libin@cs.nyu.edu}
\author[nyu]{Abtin Rahimian} \ead{arahimian@acm.org}
\author[nyu]{Denis Zorin} \ead{dzorin@cs.nyu.edu}
\address[nyu]{Courant Institute of Mathematical Sciences, New York
  University, New York, NY 10003}
\begin{abstract}
We present an efficient, accurate, and robust method for simulation of
dense suspensions of deformable and rigid particles immersed in
Stokesian fluid in two dimensions.
We use a well-established boundary integral formulation for the
problem as the foundation of our approach.
This type of formulations, with a high-order spatial discretization
and an implicit and adaptive time discretization, have been shown to
be able to handle complex interactions between particles with high
accuracy.
Yet, for dense suspensions, very small time-steps or expensive
implicit solves as well as a large number of discretization points are
required to avoid non-physical contact and intersections between
particles, leading to infinite forces and numerical instability.

Our method maintains the accuracy of previous methods at a
significantly lower cost for dense suspensions.
The key idea is to ensure interference-free configuration by
introducing explicit contact constraints into the system.
While such constraints are unnecessary in the formulation, in the
discrete form of the problem, they make it possible to eliminate
catastrophic loss of accuracy by preventing contact explicitly.

Introducing contact constraints results in a significant increase in
stable time-step size for explicit time-stepping, and a reduction
in the number of points adequate for stability.
\end{abstract}
\maketitle

\section{Introduction\label{sec:intro}}
\begin{figure}[!hbt]
  \centering
  \subfloat[$t=0$  \label{sfg:sed-snap1}]{\includegraphics[angle=-90,width=.24\linewidth]{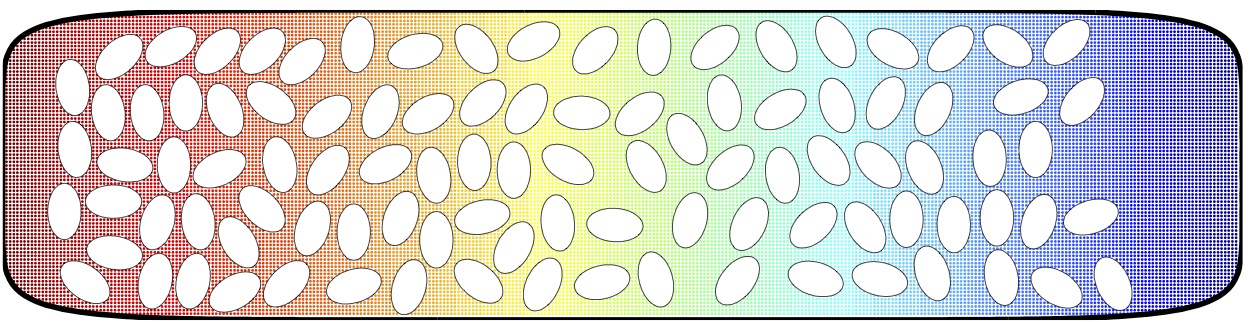}}
  \subfloat[$t=1.2$\label{sfg:sed-snap2}]{\includegraphics[angle=-90,width=.24\linewidth]{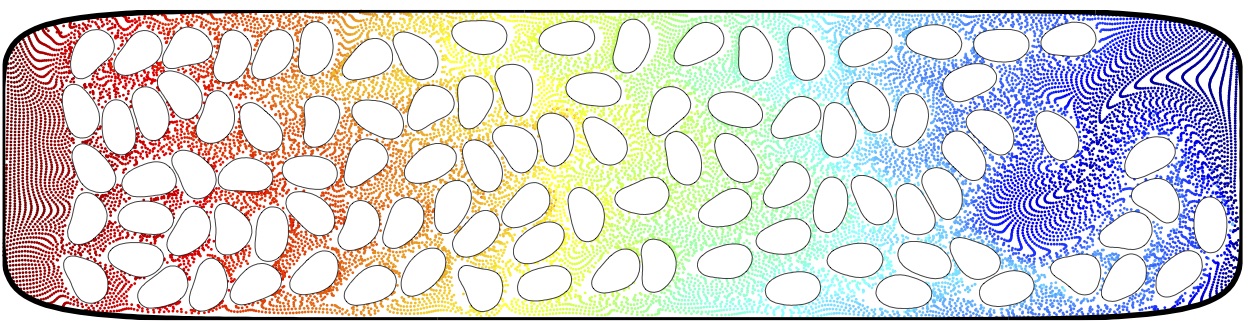}}
  \subfloat[$t=5.0$\label{sfg:sed-snap3}]{\includegraphics[angle=-90,width=.24\linewidth]{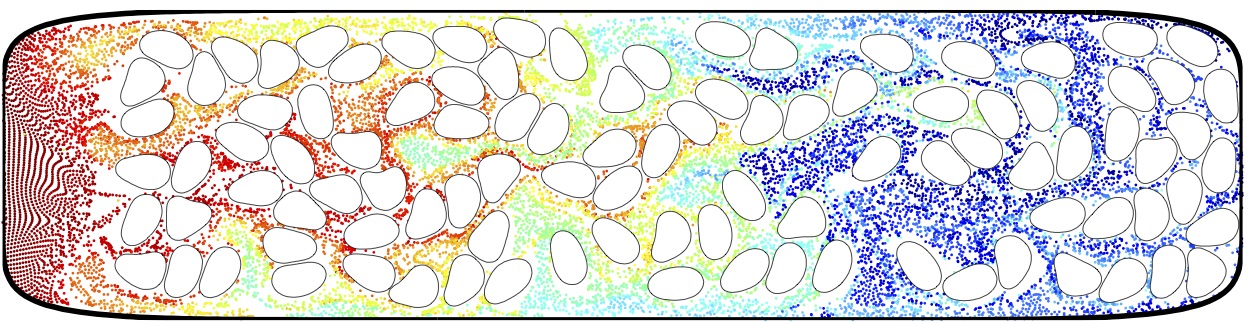}}
  \subfloat[$t=30$ \label{sfg:sed-snap4}]{\includegraphics[angle=-90,width=.24\linewidth]{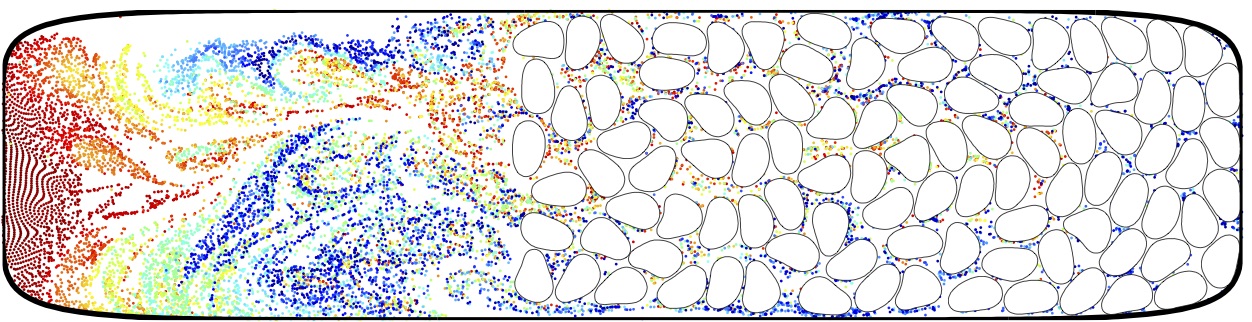}}
  \mcaption{fig:sed-snap}{Sedimentation of one hundred
    vesicles}{Vesicles are randomly placed in a container, where
    only gravitational force is present.
    \subref{sfg:sed-snap1} The initial configuration.
    The colored tracer particles are shown for visualization purposes
    and are not hydrodynamically active.
    \subref{sfg:sed-snap2} and \subref{sfg:sed-snap3}
    The intermediate states; we see that as vesicles move
    downward they induce a strong back flow.
    \subref{sfg:sed-snap4} shows the state where the simulation was
    terminated.
    As the vesicles accumulate at the bottom of the container, many
    collision areas stay active.
    Nonetheless, the vesicles are stacked stably and without any
    artifact from the collision handling.
  }
\end{figure}

Particulate Stokesian suspensions of deformable and rigid particles
are prevalent in nature and have many important industrial
applications, for example, emulsions, colloidal structures, particulate
suspensions, and blood.
Most of these examples are \emph{complex fluids}, i.e., fluids with unusual macroscopic behavior, often defying a simple constitutive-law description.
A major challenges in understanding the physics of complex fluids is
the link between microscopic and macroscopic fluid behavior.
Dynamic simulation is a powerful tool \cite{nazockdast2015b,BK2006} to
gain insight into the underlying physical principles that govern these
suspensions and to obtain relevant constitutive relationships.

Nevertheless, simulating dense suspensions of rigid and deformable
particles entails many numerical challenges, one of which is the
numerical breakdown due to collision of particles,
making robust treatment of contact vital for these simulations.
To this end, we present an efficient, accurate, and robust method for
simulation of dense suspensions in Stokesian fluid in \twod~(e.g., \cref{fig:sed-snap}),
which is extendable to \threed.
We use the boundary integral formulation to represent the flow and
impose the contact-free condition as a constraint.
This work focuses mainly on suspensions of rigid bodies and vesicles
with high volume fractions, in which multiple particles are in
contact or near-contact.

Vesicles are closed deformable membranes suspended in a viscous medium.
The dynamic deformation of vesicles and their interaction with the
Stokesian fluid play an important role in many biological phenomena.
They are used to understand the properties of
biomembranes \cite{ghigliotti2011,PhysRevLett.96.028104}, and to
simulate the motion of blood cells, in which vesicles with moderate
viscosity contrast are used to model red blood cells and high
viscosity contrast vesicles or rigid particles are used to model white
blood cells \cite{Basu20111604}.

Boundary integral formulations offers a natural approach for accurate simulation of
vesicle flows, by reducing the problem to solving equations on surfaces, and eliminating the need for discretizing changing 3D volumes.
However, in non-dilute suspensions, these methods are hindered by difficulties:
inaccuracies in computing \ns integrals, and artificial force singularities caused by
(non-physical) intersection of particles.
Contact situations in the Stokesian particulate flows occur frequently
when the volume fraction of suspensions is high, viscosity contrast of
vesicles is high, or rigid particles are present.
The dynamics of particle collision in Stokes flow are governed by the
lubrication film formation and drainage, which has a time scale much
shorter than that of the flow \cite{Frostad2013}.
Solely relying on the hydrodynamics to prevent contact requires the
accurate solution of the flow in the lubrication film, which in turn
entails very fine spatial and temporal resolution as well as expensive
implicit time-stepping\emdash/imposing excessive computational burden
as the volume fraction increases.

While adaptive time-stepping \cite{Quaife2014, Quaife2014b} goes a
long way in maintaining stability and efficiency in dilute suspensions,
the time-step is determined by the closest pair of vesicles, and tends
to be uniformly small for dense suspensions.

In this work we take a different approach: we augment the governing equations with the
contact constraint.  While from the point of view of the physics of the problem such a constraint is redundant, as non-penetration is ensured by fluid forces, in numerical context it plays an important role, improving both robustness and accuracy of simulations.
Typically, a contact \andor{law}{constraint} is characterized by
conditions of non-penetration, no-adhesion as well as a mechanical
complementarity condition, i.e., the contact force is zero when there is
no collision. These three conditions are known as Signorini conditions in the
context of contact mechanics or \abbrev{KKT} conditions in the
context of constrained optimization \cite{wriggers2006, nocedal2006}.

\subsection{Our contributions}
Contact constraints ensure that the discretized system remains
intersection-free, even for relatively coarse spatial and temporal
discretizations.
The resulting problem is a Nonlinear Complementarity Problem~(\abbrev{NCP}),
which we linearize and solve using an iterative
method that avoids explicit construction of full matrices.
We describe an implicit-explicit time-stepping scheme, adapting
Spectral Deferred Correction (\abbrev{SDC}) to our constrained
setting, \cref{sec:time-discretization}.

In our approach, the minimum distance between vesicles is controlled
by constraints, and is independent of the temporal resolution.
While this requires solving additional auxiliary equations for the
constraint forces at every step, this additional cost is more than
compensated by the ability of our method to maintain larger
time-steps, and lower spatial resolutions for a given target error.

For high volume fraction, our method makes it possible to increase the
step size by at least an order of magnitude, and the simulation remain stable
even for relatively coarse spatial discretizations (16 points per
vesicle, versus at least 64 needed for stability without contact
resolution).

\subsection{Synopsis of the method}
We use the boundary integral formulation based on
\cite{Veerapaneni2009,Rahimian2010,Quaife2014}; the basic formulation
uses integral equation form of the problem and includes the effects of
the viscosity contrast, fixed boundaries, as well as deformable and
rigid moving bodies.
We add contact constraints to this formulation, as an inequality
constraint on a gap function that is based on \emph{space-time
  intersection volume} \cite{Harmon2011}.
The contact force is then parallel to the gradient of this volume with
the Lagrange multiplier as its magnitude.

In case of multiple simultaneous contact, this leads to a Nonlinear
complementarity problem (\define{NCP}) for the Lagrange multipliers,
which we solve using a Newton-like matrix-free method, as a sequence
of Linear Complementarity Problems (\abbrev{LCP}) \cite{cottle2009,
  Erleben2013}, solved iteratively using \gmres.
The spectral Fourier bases are used for spatial discretization.
For time stepping, we use semi-implicit backward Euler or semi-implicit Spectral
Deferred Correction (\abbrev{SDC}).

\subsection{Related work}

\para{Related work on Stokesian particle flows}
Stokesian particle models are employed to theoretically and
experimentally investigate the properties of biological membranes
\cite{sackmann1996}, drug-carrying capsules \cite{sukumaran2001}, and
blood cells \cite{noguchi2005, pozrikidis1990}.
There is an extensive body of work on numerical methods for Stokesian
particulate flows and an excellent review of the literature up to 2001
can be found in \citet{pozrikidis2001a}.
Reviews of later advances can be found in
\cite{Veerapaneni2009, Rahimian2010, rahimian2015}.
Here, we briefly summarize the most important  numerical methods and
discuss the most recent developments.

Integral equation methods have been used extensively for the
simulation of Stokesian particulate flows such as droplets and bubbles
\cite{rallison1978, Zhou1993, loewenberg1998, loewenberg1997},
vesicles \cite{pozrikidis1990, freund2007, Veerapaneni2009, sohn2010,
  Rahimian2010, Farutin2013a, Zhao2011c, Zhao2013, rahimian2015}, and
rigid particles \cite{youngren1975, Power1993, Power1987}.
Other methods\emdash/such as phase-field approach \cite{biben2005,
  du2008}, immersed boundary and front tracking methods
\cite{Kim2010,Yazdani2012}, and level set method
\cite{laadhari2014}\emdash/are used by several authors for the
simulation of particulate flows.

For certain flow regimes, near interaction and collision of particles
has been a source of difficulty, which was addressed either by spatial
and temporal refinement to resolve the correct dynamics (increasing
the computational burden) or by the introduction of repulsion forces
(making the time-stepping stiff).

\citet{Sangani1994} presented a framework for dynamic simulation of
rigid particles with spherical or cylindrical shapes, in which the
lubrication forces were included directly by putting Stokes doublets
at the contact midpoint.
The magnitude of lubrication force was computed using asymptotic
analysis.
To maintain the accuracy in the interaction of deformable drops,
\cite{Zinchenko1997, loewenberg1997, Zinchenko2006} resorted to
time-step refinement where the time step is kept proportional to
particle distance~$d$.
\citet{Zinchenko1997,Zinchenko2006} keep the time step proportional
to~$\sqrt{d}$.
\citet{loewenberg1997} adjust both the grid spacing around the contact
region and the time-step to be proportional to~$d$.
\citet{freund2007} resorted to repulsion force to avoid contact in a
\twod particulate flow.
In a later work for \threed, \citeauthor{freund2007} and
coauthors~\cite{Zhao2010} observed that significantly larger repulsion
force density are needed in three dimensions, as the total repulsion force
is distributed over a smaller region, when measured as a fraction of
the total surface area/length. Consequently, they used  a purely kinematic collision handing, in which, after each time-step, the intersecting points are moved outside.

\citet{Quaife2014b} applies adaptive time-stepping and backtracking
to resolve collisions.  Similarly, \citet{Ojala2014} present an interesting
integral equation method for the flow of droplets in two dimensions with a
specialized quadrature scheme for accurate \ns evaluation enabling
simulation of flows with close to touching particles.
While methods using adaptivity both in space and time are the most robust and accurate, they incur excessive cost as means of collision handling.

\para{Related work on contact response}
A broad range of methods were developed for collision detection and
response.  While the work in contact mechanics often focuses on capturing
the physics of the contact correctly (e.g., taking into account friction effects),
the work in computer graphics literature emphasizes robustness and efficiency.
In our context, robustness and efficiency are particularly important, as we aim to
model vesicle flows with high density and large numbers of particles.
Physical correctness has a somewhat different meaning: as we know that
if the forces and surfaces in the system are resolved with high accuracy, the contacts
would not occur, our primary emphasis is on reducing the impact of the
artificial forces associated with contacts on the system.

There is an extensive literature on contact handling in
\emph{computational contact mechanics}  mainly in the context of
\abbrev{FEM} mechanical and thermal analysis
\cite{johnson1987, wriggers1995, wriggers2006, fischer2005, tur2009, puso2008, krause2002}.
\citet{wriggers1995,wriggers2006} presents in-depth reviews of the
contact mechanics framework.
The works in contact mechanics literature can be categorized
based on their ability in handling large deformations \andor{and}{or}
tangential friction.
In some of the methods, to simplify the problem, small deformation
assumption is used to predefine the active part of the boundary as
well as to align the \abbrev{FEM} mesh.
Numerical methods for contact response can be categorized as
\begin{inparaenum}[(i)]
\item penalty forces,
\item impulse/kinematic responses, and
\item constraint solvers.
\end{inparaenum}

From algorithmic viewpoint, contact mechanics methods in
\abbrev{FEM} include:
\begin{inparaenum}[(i)]
\item Node-to-node methods where the contact between nodes is only
  considered.
  The \abbrev{FEM} nodes of contacting bodies need
  to aligned and therefore this method is only applicable to small
  deformation.
\item Node-to-surface methods check the collision between predefined
  set of nodes and segments. Similar to node-to-node methods, these
  methods can only handle small deformations.
\item Surface-to-surface methods, where the contact constraint is
  imposed in weak form.
  In contrast to the two previous class of algorithms, methods in this
  class are capable of handling large deformations.
  Mortar Method is well-known within this class of algorithms
  \cite{krause2002,fischer2005,puso2008,tur2009}.
  The Mortar Method was initially developed for connecting
  different non-matching meshes in the domain decomposition approaches
  for parallel computing, e.g., \cite{puso2004}.
\end{inparaenum}

In these methods, no-penetration is either enforced as a constraint
using a Lagrangian (identified with the contact pressure) or penalty
force based on a gap function.
To the best of our knowledge, for contact mechanics problems, a signed
distance between geometric primitives is used as the gap function, in
contrast to our approach where we use space-time interference volume.

\citet{fischer2005} present a frictionless contact resolution
framework for \twod finite deformation using Mortar Method using penalty
force or Lagrange multiplier.
\citet{tur2009} use similar method for frictional contact in \twod.
\citet{puso2008} use Mortar Method for large deformation contact using
quadratic element.

Our problem has similarities to large-deformation frictionless
contact problems in contact mechanics. An important difference however, is the presence of
fluid, which plays a major role in contact response.

Application of boundary integral methods in contact mechanics is
rather limited compared to the \abbrev{FEM} methods
\cite{eck1999, gun2004}.
\citet{eck1999} used Boundary Element Method for the static contact
problem where Coulomb friction is presented.
\citet{gun2004} solved static problem with load increment and contact
constraint on displacement and traction.

In computer graphics literature, a set of commonly used and efficient methods are based on
\cite{Provot1997}, a method for the collision handling of mass-spring
cloth models.
To ensure that the system remains intersection-free, zones of impact
are introduced and rigid body motion is enforced in each zone of
impact; while this method works well in practice, its effects on the
physics of the objects are difficult to quantify.

Penalty methods are common due to the ease in their implementation,
but suffer from time-stepping stiffness \andor{and}{or}
the lack of robustness.
\citet{Baraff1998} uses implicit time-stepping coupled
with repulsion force equal to the variation of the quadratic
constraint energy with respect to control vertices.
Soft collisions are handled by the introduction of damped spring and
rigid collisions are enforced by modification to the mass matrix.
\citet{Faure2008} introduced Layered Depth Images to allow
efficient computation of the collision volumes and their gradients using \abbrev{GPU}s.
A penalty force proportional to the gradient is used to resolve collisions.
However, the stiffness of the repulsion force varies greatly (from
$\sci{5}$ to $\sci{10}$) in their experiments.
To address these difficulties, \citet{Harmon2009} present
a framework for robust simulation of contact mechanics using penalty
forces through asynchronous time-stepping, albeit at a significant
computational cost.
Alternatively, one can view collision response as an instantaneous
reaction (an impulse), i.e., an instantaneous adjustment of the velocities.
However, such adjustments are often problematic in the case of multiple contacts,
as these may lead to a cyclic ``trembling'' behavior.

Our method belongs to a large family of \emph{constraint-based} methods, which are increasingly the standard approach to contact handling in graphics. This set of methods
meets our goals of providing robustness and improving efficiency of contact response,
while minimizing the impact on the physics of the system (as virtual forces associated
with constraints do not do work).

\citet{Duriez2006} start from Signorini's law and derive the contact
force formulation.
The resulting equation is an \abbrev{LCP} that is solved by
Gauss--Seidel like iterations, sequentially resolving contacts until
reaching the contact free state.
\citet{Harmon2008} focus on robust treatment of collision without
simulation artifacts.  To enforce the no-collision constraint, this work
uses an impulse response that gives rise to an \abbrev{LCP} problem for its magnitude.
To reduce the computational cost, the \abbrev{LCP} solution (the Lagrange
multiplier) is approximated by solving a linear system.
\citet{Otaduy2009} uses a linear approximation to contact constraints
and a semi-implicit discretization, solving a mixed \abbrev{LCP}
problem at each iteration.

Our approach is directly based on \cite{Harmon2011} and is closest to
\cite{Allard2010}, in which the intersection volume and its gradient
with respect to control vertices are computed at the candidate step.
The non-collision is enforced as a constraint on this volume, which
lead to a much smaller system compared to distance formulation between
geometric primitives.
The constrained formulation leads to an \abbrev{LCP}
problem.
\citet{Harmon2011} assumes linear trajectory between edits and define
space-time interference volume and uses it as a gap function and we
use similar formulation to define the interference volume.

\subsection{Nomenclature}
In \cref{tbl:notation} we list symbols and operators frequently
used in this paper.  Throughout this paper, lower case letters refer
to scalars, and lowercase bold letters refer to vectors. Discretized
quantities are denoted by sans serif letters.
\begin{table}[!bt]
  \newcolumntype{s}{>{\centering}p{.16\linewidth}}
  \newcolumntype{d}{>{}p{.7\linewidth}}
  \setlength{\tabcolsep}{4pt}
  \centering
  \small
  \begin{minipage}[t]{.48\linewidth}
    \vspace{0pt}
    \begin{tabular}{sd}\toprule Symbol & Definition \\\midrule
      $\gamma_i$  & The boundary of the \ordinal{i} vesicle                      \\
      $\gamma$    & $\cup_i \gamma_i$                                            \\
      $\nu_i$     & The viscosity contrast                                       \\
      $\mu$       & Viscosity of the ambient fluid                               \\
      $\mu_i$     & Viscosity of the fluid inside \ordinal{i} vesicle            \\
      $\ten$      & Tension                                                      \\
      $\chi$      & Shear rate                                                   \\
      $\varpi_i$  & The domain enclosed by $\gamma_i$                            \\
      $\varpi$    & $\cup_i \varpi_i$                                            \\
      $d$         & Separation distance of particles                             \\
      $\minsep$   & Minimum separation distance                                  \\
      \Fs         & Tensile force                                                \\
      \Fb         & Bending force                                                \\
      \Fc         & Collision force                                              \\
      $h$         & Distance between two discretization points (on each surface) \\
      \bottomrule
    \end{tabular}
  \end{minipage}
  \begin{minipage}[t]{.5\linewidth}
    \vspace{0pt}
    \begin{tabular}{sd}\toprule Symbol & Definition \\\midrule
      $J$                       & Jacobian of contact volumes $V$                   \\
      $\vn$                     & Unit outward normal                               \\
      \vel                      & Velocity                                          \\
      \vbg                      & The background velocity field                     \\
      $V$                       & Contact volumes                                   \\
      $\vX$                     & Coordinate of a (Lagrangian) point on a surface   \\
      $\conv{G}$                & Stokes Single-layer operator                      \\
      $\conv{T}$                & Stokes Double-layer operator                      \\
      \hspace{0pt}\define{LCP}  & Linear Complementarity problem                    \\
      \hspace{0pt}\define{NCP}  & Nonlinear Complementarity Problem                 \\
      \hspace{0pt}\define{SDC}  & Spectral Deferred Correction                      \\
      \hspace{0pt}\define{STIV} & Space-Time Interference Volumes                   \\
      \hspace{0pt}\li           & Locally-implicit time-stepping                    \\
      \hspace{0pt}\lic          & Locally-implicit \emph{constrained} time-stepping \\
      \hspace{0pt}\gi           & Globally-implicit time-stepping                   \\
      \bottomrule
    \end{tabular}
  \end{minipage}
  \mcaption{tbl:notation}{Index of frequently used symbols, operators, and abbreviations}{}
\end{table}

\section{Formulation\label{sec:formulation}}
We consider the Stokes flow with $N_v$ vesicles and $N_p$ rigid
particles suspended in a Newtonian fluid which is either confined or
fills the free space, \cref{fig:schematic}.
In Stokesian flow, due to high viscosity \andor{and}{or} small length
scale, the ratio of inertial and viscous forces (The Reynolds number)
is small and the fluid flow can be described by the incompressible
Stokes equation
\begin{gather}
  -\mu \Delta \vu(\vx) + \Grad \scalar{p}(\vx) = \vector{F}(\vx),
  \quad\text{and}\quad \Div{\vu(\vx)} = 0 \qquad (\vx \in \Omega) \label{eq:stokes-incomp},\\
  \vector{F}(\vx) = \int_{\gamma} \vector{f}(\vX) \delta(\vx-\vX)\dline(\vX),
\end{gather}
where $\vector{f}$ is the surface density of the force exerted by the vesicle's membrane on
the fluid.
$\Omega$ denotes the fluid domain of interest with $\Gamma_0$ as its
enclosing boundary (if present) and $\mu$ denoting the viscosity of ambient
fluid.
If $\Omega$ is multiply-connected, its interior boundary consists of
$K$ smooth curves denoted by $\Gamma_1, \dots,\Gamma_K$.
The outer boundary $\Gamma_0$ encloses all the other connected
components of the domain.
The boundary of the domain is then denoted $\Gamma \defeq
\bigcup_k\Gamma_k$. We use $\vx$ to denote an Eulerian point in the fluid ($\vx
\in \Omega$) and $\vX$ a Lagrangian point on the vesicles or rigid
particles.
We let $\gamma_i$ denote the boundary of the \ordinal{i} vesicle ($i=1,
\dots, N_v$), $\varpi_i$ denote the domain enclosed by $\gamma_i$,
$\mu_i$ denote viscosity of the fluid inside that vesicle, and
 $\gamma \defeq \bigcup_{i} \gamma_i$.
\Cref{eq:stokes-incomp} is valid for $\vx\in\varpi_i$ by replacing
$\mu$ with $\mu_i$.
\begin{figure}[!bt]
  \centering
  {\hspace{-24pt}\includepgf{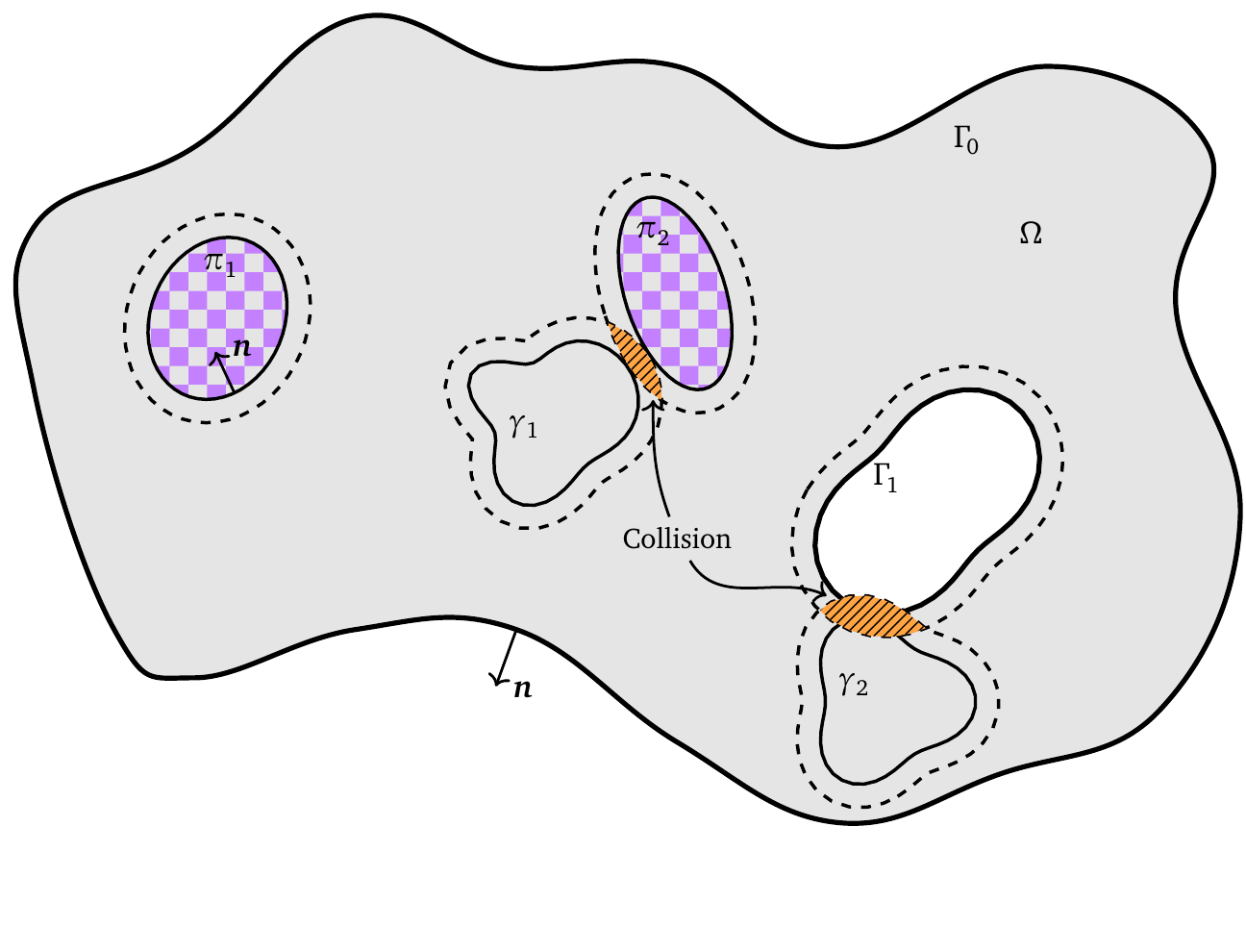}\vspace{-20pt}}
  \tikzset{external/export next=false}
  \mcaption{fig:schematic}{Schematic}{The flow domain $\Omega$
    (gray shaded area) with boundary $\Gamma_k~(k=0,\dots,K)$.
    Vesicles and rigid particles are suspended in the fluid.
    The vesicle boundaries are denoted by $\gamma_i~(i=1,\dots,N_v)$
    and the rigid particles (checkered pattern) are denoted by
    $\pi_j~(j=1,\dots,N_p)$.
    The outward normal vector to the boundaries is denotes by $\vn$.
    The dotted lines around boundaries denote the prescribed
    minimum separation distance for each of them.
    The minimum separation distance is a parameter and can be set to zero.
    In this schematic, vesicle $\gamma_1$ and particle $\pi_2$ as well
    as vesicle $\gamma_2$ and boundary $\Gamma_1$ are in contact.
    The slices of the space-time intersection volumes at the current instance are
    marked by \tikzsetnextfilename{external-schematic-pattern}\tikz \fill[preaction={fill=plt-orange},pattern=north
      east lines] (0,0) rectangle +(3ex,1.5ex);\,.
  }
\end{figure}

There are rigid particles suspended in the fluid domain.
We denote the boundary of the \ordinal{j} rigid particle  by
$\pi_j~(j=1, \dots, N_p)$ and let $\pi \defeq \bigcup_{j} \pi_j$.
The governing equations are augmented with the no-slip boundary
condition on the surface of vesicles and particles
\begin{align}
  \vu(\vX,t) &= \vX_t \qquad (\vX \in \gamma \cup \pi),
  \label{eq:interface-vel}
\end{align}
where $\vX_t \defeq \frac{\partial \vX}{\partial t}$ is the material
velocity of point $\vX$ on the surface of vesicles or particles.
The velocity on the fixed boundaries is imposed as a Dirichlet
boundary condition
\begin{equation}
  \vu(\vx) = \vector{U}(\vx) \qquad (\vx \in \Gamma).
  \label{eq:dirichlet-bc}
\end{equation}

We assume that the vesicle membrane is inextensible, i.e.,
\begin{equation}
  \label{eq:governing-membrane}
  \vX_s \cdot \vu_s = 0 \qquad  (\vX  \in \gamma),
\end{equation}
where the subscript ``$s$'' denotes differentiation with respect to the
arclength on the surface of vesicles.

Rigid particles are typically force- and torque-free.
However, surface forces may be exerted on them due to a constraint,
e.g., the contact force $\Fc$, which we will define later.
In this case, the force $\vector{F}^\pi_j$ and torque
$L^\pi_j$ exerted on the \ordinal{j} particle are the sum of
such terms induced by constraints
\begin{align}
  \vector{F}^\pi_{j} = \vector{0}, \quad \text{or} \quad
  &\vector{F}^\pi_{j} = \int_{\pi_j} \Fc(\vX) \dline \qquad (j=1,
  \dots, N_p), \label{eq:rigid-force} \\ L^\pi_{j} = \vector{0}, \quad
  \text{or} \quad &L^\pi_{j} = \int_{\pi_j} \inner{(\vX -
    \vector{c}^\pi_j)}{\perp{\Fc}(\vX)} \dline \qquad (j=1, \dots,
  N_p), \label{eq:rigid-torque}
\end{align}
where $\vector{c}^\pi_j$ is the center of mass for $\pi_j$ and for
vector $\vf=(f_1,f_2)$, $\perp{\vf} \defeq (f_2,-f_1)$.

\subsection{Contact definition\label{sec:contact-vol}}
It is known \cite{Nemer2004,Frostad2013} that the exact solution of
equations of motion,
\cref{eq:stokes-incomp,eq:interface-vel,eq:dirichlet-bc}, keeps
particles apart in finite time due to formation of lubrication film.
Thus, it is theoretically sufficient to solve the equations with an
adequate degree of accuracy to avoid any problems related to overlaps
between particles.
Nonetheless, achieving this accuracy for many types of flows (most
notably, flows with high volume fraction of particles or with complex
boundaries) may be prohibitively expensive.

With inadequate computational accuracy particles may collide with each
other or with the boundaries and depending on the numerical method
used, the consequences of this may vary.
For methods based on integral equations the consequences are
particularly dramatic, as overlapping boundaries may lead to divergent
integrals.
To address this issue, we augment the governing equations with a
contact constraint, formally written as
\begin{equation}
  \label{eq:contact-const}
  V(S,t) \le 0,
\end{equation}
where $S=\Gamma \cup \gamma \cup \pi$ denotes the boundary of the
domain and all particles.
The function $V$ is chosen in a way that $V > 0$ implies some parts
of the surface $S$ are at a distance less than a user-specified
constant \minsep.
Function $V$ may be a vector-valued function, for which the inequality
is understood component-wise.
This constraint ensures that the suspension remains contact-free
independent of numerical resolution.

For the constraint function $V,$ in addition to the basic condition
above, we choose a function that satisfies several additional
criteria:
\begin{enumerate}[(i)]
\item it introduces a relatively small number of additional
  constraints, and
\item when the function is discretized, no contacts are missed even
  for  large time steps.
\end{enumerate}
To clarify the second condition, suppose we have a small particle
rapidly moving towards a planar boundary. For a large time step, it
may move to the other side of the boundary in a single step, so any
condition that considers an instantaneous quantity depending on only
the current position is likely to miss the contact.

To this end, we use the \emph{Space-Time Interference Volumes}
(\define{STIV}) from \citet{Harmon2011} to define the function $V$ as
the area in space-time swept by the intersecting segments of boundary
over time.
To be more precise, for each point $\vX$ on the boundary $S$, consider
a trajectory $\vX(\tau)$, between a time $t_0$, for which there are no
collisions, and a time $t$.
Points $\vX(\tau)$ define a deformed boundary $S(\tau)$ for each
$\tau$.
For each point $\vX$, we define $\tau_I(\vX)$, $t_0 \leq \tau_I \leq
t$, to be the first instance for which this point comes into contact
with a different point of $S(\tau_I)$.
We let $I(\vX)$ denote the point on $S(\tau_I)$ that comes into
contact with $\vX(\tau_I)$, i.e., $\vX(\tau_I) = I(\vX)(\tau_I)$.
For points which never come into contact with other points, we set
$\tau_I = t$.
Then, the space-time interference volume for \twod problems is defined
as
\begin{equation}
  \label{eq:STV}
  V(S,t) = \int_{S(t_0)} \int_{\tau_I(\vX)}^{t}
  \left[1+\big(\inner{\vX_t(\tau)}{\vn(\vX,\tau)}\big)^2 \right]^{1/2}
\d\tau \dline
\end{equation}
where $\vector{n}(\vX,\tau)$ denotes the normal to $S(\tau)$ at
$\vX(\tau)$.  The integration is
over the initial boundary $S(t_0)$, and we use the fact that the
surface is inextensible and the surface metric does not change.
Note that, the term volume is a misnomer because \cref{eq:STV} is the
surface area swept by the intersecting segments of boundary.
An important property of this choice of function, compared to,
e.g., a space intersection volume, is that for even a very thin
object moving at high velocity, it will be proportional to the
time interval $t-t_0$.

We consider each connected component of this volume as a separate
volume, and impose an inequality constraint on each; while keeping a
single volume is in principle equivalent, using multiple volumes avoid
certain undesirable effects in discretization \cite{Harmon2011}.
Thus, $V(S,t)$ is a vector function of time-dependent dimension, with
one component per active contact region.

Depending on the context, we may omit the dependence of $V$ on $S$ and
write $V(t)$ as the contact volume function or $V(\gamma_i,t)$ for the
sub-vector of $V(S,t)$ involving surface $\gamma_i$.

In practice, it is desirable to control the minimal distance between
particles.
Therefore, we define a minimum separation distance $\minsep\ge0$ and
modify the constraint such that particles are in contact when they are
within $\minsep$ distance from each other; as shown in
\cref{fig:schematic}.
The contact volume with minimum separation distance is calculated with
the surface displaced by $\minsep$, i.e., the time $t_I$ is obtained
not from the first contact with $S(\tau)$ but rather the displaced
surface $S(\tau)+ \minsep\vn(\tau)$.
Maintaining minimum separation distance\emdash/rather than considering
pure contact only\emdash/eliminates of potentially expensive
computation of nearly singular integrals close to the surface and
improves the accuracy in semi-explicit time-stepping.

\subsection{Contact constraint\label{sec:contact-form}}
We use the Lagrange multiplier method (e.g., \cite{wriggers2006}) to
add contact constraints to the system.
While it is computationally more expensive than adding a
penalty force for the constraint (effectively, an artificial repulsion
force), it has the advantage of eliminating the need of tuning the
parameters of the penalty force to ensure that the constraint is
satisfied and keeping nonphysical forces introduced into the system to
the minimum required for maintaining the desired separation.
The constrained system can be written as
\begin{gather}
  \label{eq:opt}
  \min \int_\Omega \left(\frac{1}{2} \mu \inner{\Grad \vu}{\Grad \vu}
  - \inner{\vu}{\vector{F}} \right) \dsurf,\\ \text{subject to:}\qquad
  \begin{aligned}[t]
  \Div{\vu}(\vx) &= 0~~(\vx\in\Omega), \\ \inner{\vX_s}{\vu_s} &=
  0~~(\vX \in \gamma), \\ V(S,t) &\ge 0.\\
  \end{aligned}\notag
\end{gather}
If we omit the inequality constraint, the remaining
three equations are equivalent to the Stokes equations
\eqref{eq:stokes-incomp}. The Lagrangian for this system is
\begin{align}
  \mathcal{L}(\vu,p,\sigma,\lambda) = \int_\Omega \left(\frac{1}{2}
  \mu \inner{\Grad \vu}{\Grad \vu} - \inner{\vu}{\vector{F}} -
  p\Div{\vu} \right) \dsurf + \int_\gamma \sigma \inner{\vX_s}{\vu_s}
  \dline + \inner{\lambda}{V}.
\end{align}
The first-order optimality (\abbrevformat{KKT}) conditions yield the
following modified Stokes equation, along with the constraints
listed in \cref{eq:opt}:
\begin{gather}
  - \mu \Delta \vu + \Grad p = \vector{F}',\label{eqn:stoked-modified}\\
  \vector{F}' = \vector{F} + \int_{\gamma} \Fs \delta(\vx-\vX)\dline + \int_{S} \Fc \delta(\vx-\vX)\dline,  \label{eq:vol-force}\\
  \Fs = -(\sigma \vX_s)_s,\\
  \Fc =  \Grad V^T \lambda,\\
  \lambda \ge 0,\\
  \inner{\lambda}{V} = 0,
\end{gather}
where the last condition is the complementarity
condition\emdash/either an equality constraint is active $(V_i=0)$ or
its Lagrange multiplier is zero.
As we will see in the next section and based on \cref{eq:vol-force},
the collision force $\Fc$ is added to the traction jump across the
vesicle's interface.
For rigid particles, the contact force induces force and torque on
each particle\emdash/as given in
\cref{eq:rigid-force,eq:rigid-torque}.

It is customary to combine $V\ge 0, \lambda\ge 0$, and
$\inner{\lambda}{V} = 0$,  into one expression and write
\begin{equation}
  0 \le V(t) \quad \bot \quad \lambda \ge 0, \label{eq:signorini}
\end{equation}
where "$\bot$" denotes the complementarity condition.
These ensure that the Signorini conditions introduced in
\cref{sec:intro} are respected: contacts do not produce attraction
force ($\lambda \ge 0$) and the constraint is active ($\lambda$
nonzero) if and only if $V(t)$ is zero.
Furthermore, it follows from $V(t) \cdot \lambda = 0$ that the contact
force $\Fc$ is perpendicular to the velocity and therefore it respects
the principal of virtual work and does not add to or remove from the
system's energy.

\subsection{Boundary integral formulation}
Following the standard approach of potential theory \cite{Poz92,
  Power1993}, one can express the solution of the Stokes
boundary value problem, \cref{eqn:stoked-modified}, as a system of
singular integro-differential equations on all immersed and bounding
surfaces.

The Stokeslet tensor $\linop{G}$, the Stresslet tensor $\linop{T}$,
and the Rotlet $\linop{R}$ are the fundamental solutions of the Stokes
equation given by
\begin{align}
  \linop{G}(\vector{r}) &= \frac{1}{4\pi\mu}
  \left(-\log\norm{\vector{r}} \vector{I} + \frac{\vector{r} \otimes
    \vector{r}}{\norm{\vector{r}}^2}
  \right), \label{eq:single-layer}\\ \linop{T}(\vector{r}) &=
  \frac{1}{\pi\mu} \frac{\vector{r} \otimes \vector{r} \otimes
    \vector{r}}{ \norm{ \vector{r}}^4
  }, \label{eq:double-layer}\\ \linop{R}(\vector{r}) &= \frac{1}{4
    \pi \mu} \frac{ \perp{\vector{r}}}{\norm{\vector{r}}^2}.
\end{align}

The solution of \cref{eqn:stoked-modified} can be expressed by the
combination of single- and double-layer integrals.
We denote the \slyr integral on the vesicle surface $\gamma_i$ by
\begin{equation}
  \label{eq:SLP}
  \conv{G}_{\gamma_i}[\vector{f}](\vx) \defeq \int_{\gamma_i}
  \linop{G}(\vx - \vY)\cdot\vector{f}(\vY) \dline,
\end{equation}
where $\vector{f}$ is an appropriately defined density.
The \dlyr integral on a surface $S$ (a vesicle, a rigid
particle, or a fixed boundary) is
\begin{equation}
  \label{eq:DLP}
  \conv{T}_S[\vector{q}](\vx) \defeq \int_{S}
  \vector{n}(\vY) \cdot \linop{T}(\vector{x} -
  \vY)\cdot\vector{q}(\vY) \dline,
\end{equation}
where $\vn$ denotes the outward normal (as shown in
\cref{fig:schematic}) to the surface $S$, and $\vector{q}$ is
appropriately defined density.
When the evaluation point $\vx$ is on the integration surface,
\cref{eq:SLP} is a singular integral, and \cref{eq:DLP} is
interpreted in the principal value sense.

Due to the linearity of the Stokes equations, as formulated in
\cite{Rahimian2010, Quaife2014}, the velocity at a point
$\vx\in\Omega$ can be expressed as the superposition of velocities due
to vesicles, rigid particles, and fixed boundaries
\begin{equation}
  \label{eq:vel-super}
  \alpha \vu(\vx) = \vbg(\vx) + \vu^{\gamma}(\vx) + \vu^{\pi}(\vx) +
  \vu^{\Gamma}(\vx), \quad \vx \in \Omega, \quad \alpha
  = \begin{cases} 1 & \vx \in \Omega\backslash\gamma, \\ \nu_i & \vx
    \in \varpi_i, \\ (1+\nu_i)/2 & \vx \in \gamma_i, \end{cases}
\end{equation}
where $\vbg(\vx)$ represent the background velocity field (for
unbounded flows) and $\nu_i=\mu_i/\mu$ denotes the viscosity contrast
of the \ordinal{i} vesicle.
The velocity contributions from vesicles, rigid particles, and fixed
boundaries each can be further decomposed into the contribution of
individual components
\begin{align}
  \vu^{\gamma}(\vx) &= \sum_{i=1}^{N_v} \vu_{i}^{\gamma}(\vx), &
  \vu^{\pi}(\vx)    &= \sum_{j=1}^{N_p} \vu_{j}^{\pi}(\vx),    &
  \vu^{\Gamma}(\vx) &= \sum_{k=0}^{K} \vu_{k}^{\Gamma}(\vx).
\end{align}

To simplify the representation, we introduce the complementary
velocity for each component. For the \ordinal{i} vesicle, it is defined as
$\comp{\vu}^\gamma_i = \alpha\vu - \vu^\gamma_i$.
The complementary velocity is defined in a similar fashion for rigid
particles as well as components of the fixed boundary.

\subsubsection{The contribution from vesicles}
The velocity induced by the \ordinal{i} vesicle is expressed as an
integral \cite{Poz92}:
\begin{equation}
  \vu_{i}^{\gamma}(\vx) = \conv{G}_{\gamma_i}[\vf](\vx) +
  (1-\nu_i)\conv{T}_{\gamma_i}[\vector{u}](\vx) \qquad(\vx\in\Omega),
  \label{eq:vesicle-vel}
\end{equation}
where the \dlyr density $\vu$ is the total velocity and $\vector{f}$
is the traction jump across the vesicle membrane. Based on
\cref{eq:vol-force}, the traction jump is equal to the sum of bending,
tensile, and collision forces
\begin{equation}
 \vector{f}(\vX) = \Fb + \Fs + \Fc = -\kappa_b\vX_{ssss} - (\sigma
 \vX_s)_s + \Grad V^T \lambda\qquad(\vX \in \gamma),
\label{eq:traction-jump}
\end{equation}
where $\kappa_b$ is the membrane's bending modulus.
The tensile force $\Fs = (\sigma \vX_s)_s$ is determined by the
local inextensibility constraint, \cref{eq:governing-membrane}, and
tension $\sigma$ is its Lagrangian multiplier.

Note that \cref{eq:vesicle-vel} is the contribution from each vesicle
to the velocity field.
To obtain an equation for the interfacial velocity,
\cref{eq:vesicle-vel} is to be substituted into \cref{eq:vel-super}:
\begin{gather}
  \frac{(1+\nu_i)}{2} \vu(\vX) = \comp{\vu}_i^\gamma(\vX) +
  \conv{G}_{\gamma_i}[\vf](\vX) +
  (1-\nu_i)\conv{T}_{\gamma_i}[\vector{u}](\vX)
  \qquad(\vX\in\gamma_i), \label{eq:vesicle-bie}\\ \intertext{subject
    to local inextensibility} \vX_s \cdot \vu_s = 0 \qquad (\vX \in
  \gamma_i)\label{eq:vesicle-tesion}.
\end{gather}

\subsubsection{The contribution from the fixed boundaries\label{ssc:bdry}}
The velocity contribution from the fixed boundary can be expressed as
a \dlyr integral \cite{Power1993} along $\Gamma$.
The contribution of the outer boundary $\Gamma_0$ is
\begin{equation}
  \vu_{0}^{\Gamma}(\vx) = \conv{T}_{\Gamma_0}[\vector{\eta_{0}}](\vx)
  \qquad(\vx\in\Omega), \label{eq:bound-out}
\end{equation}
where $\vector{\eta}_0$ is the density to be determined based on
boundary conditions.
Substituting \cref{eq:bound-out} into \cref{eq:vel-super} and taking
its limit to a point on $\Gamma_0$ and using the Dirichlet boundary
condition, \cref{eq:dirichlet-bc}, we obtain a Fredholm integral
equations for the density $\vector{\eta}_0$
\begin{equation*}
    \vector{U}(\vx) - \comp{\vu}_0^\Gamma(\vx) =
    -\frac{1}{2}\vector{\eta}_0(\vx) +
    \conv{T}_{\Gamma_0}[\vector{\eta_{0}}](\vx)
    \qquad(\vx\in\Gamma_0).
\end{equation*}
However, this equation is rank deficient.  To render it invertible,
the equation is modified following \cite{Poz92}:
\begin{equation}
    \vector{U}(\vx) - \comp{\vu}_0^\Gamma(\vx) =
    -\frac{1}{2}\vector{\eta}_0(\vx) +
    \conv{T}_{\Gamma_0}[\vector{\eta_{0}}](\vx) +
    \conv{N}_{\Gamma_0}[\vector{\eta_{0}}](\vx)
    \qquad(\vx\in\Gamma_0), \label{eq:bie-outer-comp}
\end{equation}
where the operator $\conv{N}_{\Gamma_0}$ is defined as
\begin{equation}
  \conv{N}_{\Gamma_0}[\vector{\eta_{0}}](\vx) = \int_{\Gamma_0}
  \contract{[\vn(\vx)\otimes\vn(\vy)]}{\vector{\eta_{0}}(\vy)} \dline.
\end{equation}

For the enclosed boundary components $\Gamma_k~(k>0)$ to eliminate the
double-layer nullspace we need to include additional Stokeslet and Rotlet terms
\begin{equation}
  \vu_{k}^{\Gamma}(\vx) = \conv{T}_{\Gamma_k}[\vector{\eta_{k}}](\vx)
  + \linop{G}(\vx - \vector{c}_{k}^{\Gamma}) \cdot
  \vector{F}_{k}^{\Gamma} + \linop{R}(\vx - \vector{c}_{k}^{\Gamma})
  L_{k}^{\Gamma}, \qquad~(k=1, \dots, K; \vx\in\Omega),
  \label{eq:bound-in}
\end{equation}
where $\vector{c}^\Gamma_k$ is a point enclosed by $\Gamma_k$,
$\vector{F}_k^\Gamma$ is the force exerted on $\Gamma_k$, and
$L_k^\Gamma$ is the torque:
\begin{align}
  \vector{F}_{k}^{\Gamma} &= \frac{1}{|\Gamma_k|}\int_{\Gamma_k}
  \vector{\eta_k} \dline, & L_{k}^{\Gamma} &=
  \frac{1}{|\Gamma_k|}\int_{\Gamma_k} \inner{(\vX -
  \vector{c}_{k}^{\Gamma})}{\perp{\vector{\eta}}}_k \dline,
  \label{eq:bie-comp-constraint}
\end{align}
where $|\Gamma_k|$ denotes the perimeter of $\Gamma_k$.
Taking the limit at points of the surface, leads to the following integral equation:
\begin{align}
 \vector{U}(\vx) - \comp{\vu}_k^\Gamma(\vx) &=
 -\frac{1}{2}\vector{\eta}_k(\vx) +
 \conv{T}_{\Gamma_k}[\vector{\eta_{k}}](\vx) + \linop{G}(\vx -
 \vector{c}_{k}^{\Gamma}) \cdot \vector{F}_{k}^{\Gamma} +
 \linop{R}(\vx - \vector{c}_{k}^{\Gamma})
 L_{k}^{\Gamma}\qquad(\vx\in\Gamma_k). \label{eq:bie-inner-comp}
\end{align}
%
%
\Cref{eq:bie-inner-comp,eq:bie-comp-constraint} are a complete system
for \dlyr densities $\vector{\eta}_k$, forces $\vector{F}^\Gamma_k$,
and torques $L_k^\Gamma$ on each enclosed surface $\Gamma_k$.

\subsubsection{The contribution from rigid particles\label{ssc:rigid}}
The formulation for rigid particles is very similar to that of fixed
boundaries, except the force and torque are known (cf. \cref{eq:rigid-force,eq:rigid-torque}).
The velocity contribution from the \ordinal{j} rigid particle is
\begin{equation}
  \label{eq:rigid-DLP}
  \vu_{j}^{\pi}(\vx) = \conv{T}_{\pi_j}[\vector{\zeta}_{j}](\vx) +
  \vector{G}(\vx - \vector{c}_{j}^{\pi}) \cdot \vector{F}_{j}^{\pi} +
  \vector{R}(\vx - \vector{c}_{j}^{\pi}) L_{j}^{\pi},
\end{equation}
Where $ \vector{F}_{j}^{\pi}, L_{j}^{\pi}$ are, respectively, the \emph{known}
 net force and torque exerted on the particle and  $\vector{\zeta}_j$ is the unknown density.

Let $\vector{U}_j^\pi$ and $\omega_j^\pi$ be the translational
and angular velocities of the \ordinal{j} particle; then we obtain the
following integral equation for the density~$\vector{\zeta}_j$ from the limit of
\eqref{eq:rigid-DLP}:
\begin{equation}
  \label{eq:rigid-fredholm}
  \vector{U}_j^\pi + \omega_j^\pi \perp{ (\vX-\vector{c}_j^\pi)}
  - \comp{\vu}^\pi_j(\vX) = - \frac{1}{2}\vector{\zeta}_j(\vX) +
  \conv{T}_{\pi_j}[\vector{\zeta}_{j}](\vX) + \vector{G}(\vX -
  \vector{c}_{j}^{\pi}) \cdot \vector{F}_{j}^{\pi} + \vector{R}(\vX -
  \vector{c}_{j}^{\pi}) L_{j}^{\pi}.
\end{equation}
where
\begin{align}
  \vector{F}_j^{\pi} &= \frac{1}{|\pi_j|}\int_{\pi_j}\vector{\zeta}_j
  \dline , & L_{j}^{\pi} &= \frac{1}{|\pi_j|}\int_{\pi_j}
  \inner{(\vector{Y} - \vector{c}_j^{\pi})}{\perp{\vector{\zeta}}_j}
  \dline
  \label{eq:rigid-const}
\end{align}
where $\vector{c}_{j}^{\pi}$ is the center of $\ordinal{j}$ rigid
particle.  \Cref{eq:rigid-fredholm,eq:rigid-const} are used to solve
for the unknown densities $\vector{\zeta}_j$ as well as the unknown
translational and angular velocities of each particle. Note that the
objective of \cref{eq:bie-comp-constraint,eq:rigid-const} is remove
the null space of the \dlyr and therefore the left-hand-side (i.e. the
image of the null function) can be chosen rather arbitrarily.

\subsection{Formulation summary}
The formulae outlined above govern the evolution of the suspension.
The flow constituents are hydrodynamically coupled through the
complementary velocity (i.e., the velocity from all other
constituents).
Given the configuration of the suspension, the
unknowns are:
\begin{itemize}
  \item Velocity $\vu(\vX)$ and tension $\sigma$ of vesicles'
    interface determined by \cref{eq:traction-jump,eq:vesicle-bie,eq:vesicle-tesion}.
    The velocity is integrated for the vesicles' trajectory using
    \cref{eq:interface-vel}.
    %
  \item The \dlyr density on the enclosing boundary $\vector{\eta}_0$
    as well as the \dlyr density $\vector{\eta}_k~(k=1,\dots,K)$,
    force $\vector{F}_k^\Gamma$, and torque $L_k^\Gamma$ on the
    interior boundaries determined by
    \cref{eq:bie-outer-comp,eq:bie-inner-comp,eq:bie-comp-constraint}.
    Note that the collision constraint does not enter the formulation
    for the fixed boundaries and when a particle collides with a fixed
    boundary, the collision force is only applied to the particle.
    The unknown force and torque above can be interpreted as the
    required force to keep the boundary in place.
  \item Translational $\vector{U}_j^\pi$ and angular $\omega_j^\pi$
    velocities of rigid particles $(j=1,\dots,N_p)$ as well as \dlyr
    densities $\vector{\zeta}_j$ on their boundary determined by
    \cref{eq:rigid-fredholm,eq:rigid-const}.
    Where the force and torque are either zero or determined by the
    collision constraint \cref{eq:rigid-force,eq:rigid-torque}.
\end{itemize}

This system is constrained by Signorini (\abbrev{KKT}) conditions for
the contact, \cref{eq:signorini},
which is used to compute $\lambda$, the strength of the contact force.

In the referenced equations above, the complementary velocity is
combination of velocities given in
\cref{eq:vesicle-vel,eq:bound-out,eq:bound-in,eq:rigid-DLP}.

\section{Discretization and Numerical Methods\label{sec:numerical}}
In this section, we describe the numerical algorithms required for
solving the dynamics of a particulate Stokesian suspension.
We use the spatial representation and integral schemes in \cite{Rahimian2010}.
We also adapt the spectral deferred correction time-stepping in
\cite{Quaife2014a,Quaife2014b} to the local implicit time-stepping schemes.
Furthermore, we use piecewise-linear discretization of curves to calculate
the space-time contact volume $V(\gamma,t)$, \cref{eq:STV}, as in \cite{Harmon2011}.
To solve the complementarity problem resulting from the contact
constraint, we use the minimum-map Newton method discussed in
\cite{Erleben2013}.

The key difference, compared to previous works is that at every time
step instead of solving a linear system we solve a nonlinear
complementarity problem (\abbrev{NCP}).
The \abbrev{NCP}s are solved iteratively, using a Linear
Complementarity Problem (\abbrev{LCP}) solver.
We refer to these iterations as contact-resolving iterations, in
contrast to the outer time-stepping iterations.

For simplicity, we describe the scheme for a system including vesicles
only, without boundaries or rigid particles.  Adding these requires
straightforward modifications to the equations.  In the following
sections, we will first summarize the spatial discretization, then
discuss the \abbrev{LCP} solver, and close with the time
discretization with contact constraint.

\subsection{Spatial discretization}
All interfaces are discretized with $N$ uniformly-spaced discretization points
\cite{Rahimian2010}. The number of points on each curve is typically
different but for the sake of clarity we use $N$.
The distance between discretization points does not change with time
over the curves due to rigidity of particles or the local
inextensibility constraint for vesicles.
Let $\vX(s)$, with $s \in (0,L]$, be a parametrization of the
  interface $\gamma_i~(\text{or}~\pi_j)$, and let $\{s_k =
  kL/N\}_{k=1}^{N}$ be $N$ equally spaced points in arclength
  parameter, and $\vXd_k \defeq \vX(s_k)$ the corresponding
  material points.

\para{High-order discretization for force computation} We use the
Fourier basis to interpolate the positions and forces associated with
sample points, and \abbrev{FFT} to calculate the derivatives of all
orders to spectral accuracy.
We use the hybrid Gauss-trapezoidal quadrature rules of
\cite{Alpert99} to integrate the \slyr potential for
$\vX\in\gamma_i$
\begin{equation}
  \conv{G}_{\gamma_i}[\vector{f}](\vX) \approx
  \linopd{G}_{\gamma_i}[\vectord{f}](\vXd) \defeq \sum_{\ell=1}^{N+M}
  w_\ell \contract{\linopd{G}(\vXd - \vYd_\ell)}{\vectord{f}(
    \vectord{Y}_\ell)},
\end{equation}
where $w_\ell$ are the quadrature weights given in \cite[Table
  8]{Alpert99} and $\vYd_\ell$ are quadrature points.
For $\vX\in\gamma_i$, the linear operator
$\linopd{G}_{\gamma_i}$ is a matrix,
we denote the \slyr potential on $\gamma_i$ as
$\linopd{G}_{\gamma_i}\vectord{f}$.

The \dlyr kernel $\contract{\vn(\vY)}{\linop{T}(\vX - \vY)}$ in
  \cref{eq:DLP} is non-singular in two dimensions
\[\lim_{\gamma_i\ni\vX\to\vY} \vn(\vY) \cdot \linop{T}(\vX - \vY) = -\frac{\kappa}{2\pi\mu}\vector{t}\otimes\vector{t},\]
where $\vector{t}$ denotes the tangent vector.
Therefore, a simple uniform-weight composite trapezoidal quadrature
rule has spectral accuracy in this case.
Similar to the \slyr case, we denote the discrete \dlyr potential
on $\gamma_i$ by $\linopd{T}_{\gamma_i}\vectord{u}$.
We use the nearly-singular integration scheme described in
\cite{Quaife2014} to maintain high integration accuracy for particles
closely approaching each other.


\para{Piecewise-linear discretization for constraints.}
While the spectral spatial discretization is used for most
computations, it poses a problem for the minimal-separation constraint
discretization.
Computing parametric curve intersections, an essential step in the
\abbrev{STIV} computation, is relatively expensive and difficult to
implement robustly, as this requires solving nonlinear equations.
We observe that the impact of the separation distance on the overall
accuracy is low in most situations, as explored in \cref{sec:results}.
Thus, rather than enforcing the constraint as precisely as allowed by
the spectral discretization, we opt for a low-order, piecewise-linear
discretization in this case, and use an algorithm that ensures that
\emph{at least} the target minimal separation is maintained, but may
enforce a higher separation distance.

For the purpose of computing \abbrev{STIV} and its gradient, we use
$L(\vXd,r)$, the piecewise-linear interpolant of $r$ times refinement
of points\emdash/the upsampled points correspond to arclength values
with spacing $L/(N2^r)$, with $r$ determined dynamically.

For discretized computations, we set the separation distance to
$(1+2\alpha)\minsep$, where $\minsep$ is the target minimum separation
distance.
We choose $r$ such that $\norm[\infty]{L(\vXd,r)-\vX} < \alpha
\minsep$.
Our \abbrev{NCP} solver, described below, ensures that the separation
between $L(\vXd,r)$ is $(1+2\alpha) \minsep$ at the end of a single
time step. We choose $\alpha = 0.1$ which requires $r=1$ in our
experiments; smaller values of $\alpha$ require more refinement, but
enforce the constraint more accurately.

At the end of the time step, the minimal-separation constraint ensures
that $L(\vXd,r)(s)$, for any $s$, is at least at the distance
$(1+2\alpha) \minsep$ from a possible intersection if its trajectory is
extrapolated linearly.  By computing the upper bounds on the
difference between the $\vX(s)$ and $L(\vXd,r)(s)$ at the beginning of
the time step, and interpolated velocities, we obtain a lower bound on
the actual separation distance $d'$ for the spectral surface $\vX(s)$.
If $d' < \minsep$, we increase $r$, and repeat the time step.  As the
piecewise linear approximation converges to the spectral boundary
$\vX$, and so do the interpolated velocities. In practice, we have not
observed a need for refinement for our choice of $\alpha$.

For the piecewise-linear discretization of curves,
the space-time contact volume $V(\gamma,t)$, \cref{eq:STV}, and its
gradient are calculated using the definitions and algorithms in \cite{Harmon2011}.
Given a contact-free configuration and a candidate configuration for the next time step,
we calculate the discretized space-time contact volume as the sum of edge-vertex
contact volumes $V = \sum_k V_k(\vectord{e},\vXd)$.
We use a regular grid of size proportional to the average boundary
spacing to quickly find potential collisions.
For all vertices and edges, the bounding box enclosing their initial
(collision-free) and final (candidate position) location is formed and
all the grid boxes intersecting that box are marked.
When the minimal separation distance $\minsep>0$, the bounding box is enlarged by $\minsep$.
For each edge-vertex pair $\vectord{e}(\vXd_i,\vXd_{i+1})$ and
$\vXd_j$, we solve a quartic equation to find their earliest contact
time $\tau_I$ assuming linear trajectory between initial and candidate
positions.
We calculate the edge-vertex contact volume using \cref{eq:STV}:
\begin{equation}
  V_k(\vectord{e},\vXd) = (t-\tau_I)(1+(\vXd(\tau_I)\cdot\vectord{n}(\tau_I))^2)^{1/2} |\vector{e}|,
\end{equation}
where $\vectord{n}(\tau_I)$ is the normal to the edge $\vectord{e}(\tau_I)$.
For each edge-vertex contact volume, we calculate the gradient respect to the vertices
$\vXd_i$, $\vXd_{i+1}$ and $\vXd_{j}$, summing over all the edge-vertex contact pairs we
get the total space-time contact volume and gradient.

\subsection{Temporal discretization\label{sec:time-discretization}}
Our temporal discretization is based on  the locally-implicit
time-stepping scheme in \cite{Rahimian2010}\emdash/ adapting the
Implicit--Explicit (\abbrev{IMEX}) scheme \cite{Ascher1995a} for
interfacial flow\emdash/ in which we treat intra-particle interactions
implicitly and inter-particle interactions explicitly.
We combine this method with the minimal-separation constraint.
We refer to this scheme as \emph{constrained locally-implicit} (\lic)
scheme.
For comparison purposes, we also consider the same scheme without
constraints (\li) and the \emph{globally semi-implicit} (\gi) scheme,
where all interactions treated implicitly \cite{rahimian2015}.
From the perspective of boundary integral formulation, the
distinguishing factor between \li and \lic is the extra traction jump
term due to collision.
Schemes \andor{\li}{\lic} and \gi differ in their explicit or implicit
treatment of the complementary velocities.

While treating the inter-vesicle interactions explicitly may result in
more frequent violations of minimal-separation constraint, we
demonstrate that in essentially all cases the \lic scheme is
significantly more efficient than both the \gi and \li schemes because
these schemes are costlier and require higher spatial and temporal resolution to
prevent collisions.

We consider two versions of the \lic scheme, a simple first-order
Euler scheme and a spectral deferred correction version.
A first-order backward Euler \lic time stepping formulation for
\cref{eq:vesicle-bie} is
\begin{align}
  \frac{1+\nu_i}{2} \vectord{u}_{i}^{+} &= \comp{\vectord{u}}_i^\gamma
  + \linopd{G}_{\gamma_i}\vectord{f}_{i}(\vXd_{i}^{+}, \sigma_i^{+},
  \lambda^{+}) + (1-\nu_i)\linopd{T}_{\gamma_i}\vectord{u}_{i}^{+},
  \label{eq:vesicle-bie-disc}\\
  \inner{\vXd_{i,s}}{\vectord{u}_{i,s}^{+}} &= 0 ,\\
  \vectord{f}_i(\vXd_{i}^{+}, \sigma_i^+, \lambda^+) &=
  -\kappa_b\vXd_{i,ssss}^{+} - (\sigma^{+} \vXd_{i,s})_s +
  ((\Grad \discrete{V}^{+})^T \lambda^+ )_i,\\
  0 \le\discrete{V}(\gamma;t^{+}) \quad & \bot \quad \lambda^{+} \ge 0,
\end{align}
where the implicit unknowns to be solved for at the current step are
marked with superscript ``$^{+}$''.
The position and velocity of the points of $\ordinal{i}$ vesicle are
denoted by $\vXd_i$, $\vectord{u}_i^+ = (\vXd_i^+ - \vXd_i)/\dt$, and
$\vectord{f}_i$ is the traction jump on the $\ordinal{i}$ vesicle
boundary.
$\discrete{V}(\gamma;t^{+})$ is the \abbrev{STIV} function.

\subsubsection{Spectral Deferred Correction}
We use spectral deferred correction~(\abbrev{SDC}) method
\cite{Quaife2014a,Quaife2014b}  to get a better
stability behavior compared to the basic backward Euler scheme
described above.
We use \abbrev{SDC} both for \li and \lic time-stepping.
To obtain the \abbrev{SDC} time-stepping equations, we reformulate
\cref{eq:interface-vel} as a Picard integral
\begin{equation}
\label{eq:picard}
  \vX(t_{n+1}) = \vX(t_n) + \int_{t_n}^{t_{n+1}} \vector{u}\left(\vX(\tau),
  \sigma(\tau),\lambda(\tau)\right) \d \tau,
\end{equation}
where the velocity satisfies \cref{eq:vesicle-bie,eq:vesicle-bie-disc}.
In the \abbrev{SDC} method,  the integral in \cref{eq:picard} is first discretized
with $p+1$ Gauss-Lobatto quadrature points.
Each iteration starts with
$p$ \emph{provisional positions} $\widetilde{\vXd}$ corresponding to
times $\tau_i$ in the interval $[t_n,t_{n+1}]$;  $t_n = \tau_0 < \cdots < \tau_p = t_{n+1}$.
 Provisional tensions $\widetilde{\scalard{\sigma}}$
 and provisional $\widetilde{\scalard{\lambda}}$ are defined similarly.
The \abbrev{SDC} method iteratively corrects the provisional positions $\widetilde{\vXd}$
with the error term $\widetilde{\vectord{e}}$, which is solved using the residual $\widetilde{\vectord{r}}$ resulting from the provisional solution as defined below. The residual is given by:
\begin{equation}
  \label{eq:sdc-res}
  \widetilde{\vector{r}}(\tau) = \widetilde{\vX}(t_n)
  - \widetilde{\vX}(\tau) + \int_{t_n}^{\tau}  \widetilde{\vector{u}}(\theta)\d \theta.
\end{equation}

After discretization, we use $\widetilde{\vXd}^{w,m}$,
to denote the provisional position at \ordinal{m} Gauss-Lobatto point
after $w$ \abbrev{SDC} passes.
The error term $\widetilde{\vectord{e}}^{w,m}$
denotes the computed correction to obtain \ordinal{m} provisional
position in \ordinal{w} iteration.
The \abbrev{SDC} correction iteration is defined by
\begin{align}
  \widetilde{\vXd}^{w,m} &= \widetilde{\vXd}^{w-1,m} + \widetilde{\vectord{e}}^{w,m},\notag\\
  \widetilde{\sigma}^{w,m} &= \widetilde{\scalard{\sigma}}^{w-1,m} +
  \widetilde{\scalard{e}}^{w,m}_{\sigma}, \label{eq:sdc-correct}\\
  \widetilde{\lambda}^{w,m} &= \widetilde{\scalard{\lambda}}^{w-1,m} +
  \widetilde{\scalard{e}}^{w,m}_{\lambda}\notag.
\end{align}
Setting $\widetilde{\vXd}^{0,m}$ to zero, the first \abbrev{SDC} pass
is just backward Euler time stepping to obtain nontrivial provisional solutions.
Beginning from the second pass, we solve for the error term as corrections.

Denote $(\alpha \linopd{I} - (1-\nu)\linopd{T}_{\widetilde{\gamma}^{w-1,m}})$ by
$\linopd{D}_{\widetilde{\gamma}^{w-1,m}}$. Following \cite{Quaife2014a,Quaife2014b},
we solve the following equation for the error term:
\begin{equation}
  \label{eq:sdc-error-disc}
  \alpha \frac{\widetilde{\vectord{e}}^{w,m}
  - \widetilde{\vectord{e}}^{w,m-1}}{\Delta \tau} =
  \linopd{D}_{\widetilde{\gamma}^{w-1,m}}
  \left(\frac{\widetilde{\vectord{r}}^{w-1,m}
  - \widetilde{\vectord{r}}^{w-1,m-1}}{\Delta \tau}\right)
  + \linopd{G}_{\widetilde{\gamma}^{w-1,m}}
  \vectord{f}(\widetilde{\vectord{e}}^{w,m},\widetilde{\scalard{e}}^{w,m}_{\sigma},
  \widetilde{\scalard{e}}^{w,m}_{\lambda})
  + (1-\nu)\linopd{T}_{\widetilde{\gamma}^{w-1,m}}
  \left(\frac{\widetilde{\vectord{e}}^{w,m}
  - \widetilde{\vectord{e}}^{w,m-1}}{\Delta \tau}\right).
\end{equation}
\cref{eq:sdc-error-disc} is the identical to  \cref{eq:vesicle-bie-disc},
except the right-hand-side for \cref{eq:sdc-error-disc} is obtained from the residual while
the right-hand-side for \cref{eq:vesicle-bie-disc} is the complementary velocity.
The residual $\widetilde{\vectord{r}}^{w,m}$ is obtained using
a discretization of  \cref{eq:sdc-res}:
\begin{equation}
  \label{eq:sdc-res-disc}
  \widetilde{\vectord{r}}^{w,m} = \widetilde{\vXd}^{w,0}
  - \widetilde{\vXd}^{w,m} + \sum_{l = 0}^{p} w_{l,m} \widetilde{\vectord{u}}^{w,l}.
\end{equation}
where $w_{l,m}$ are the quadrature weights for Gauss-Lobatto points, whose
quadrature error is $\mathcal{O}(\Delta t^{2p-3})$.
In addition to the \abbrev{SDC} iteration, \cref{eq:sdc-error-disc},
we also enforce the inextensibility constraint
\begin{equation}
  \label{eq:sdc-error-sigma-disc}
  \widetilde{\vXd}^{w-1,m}_s \cdot \widetilde{\vectord{u}}^{w,m}_s = 0,
\end{equation}
and the contact complementarity
\begin{equation}
  \label{eq:sdc-error-lambda-disc}
  0 \le \discrete{V}(\widetilde{\gamma}^{w,m}) \quad  \bot \quad
  \widetilde{\scalard{e}}_\lambda^{w,m} \ge 0.
\end{equation}

In evaluating the residuals using \cref{eq:sdc-res-disc}, provisional
velocities are required. In the \gi scheme \cite{Quaife2014a}, all the
interactions are treated implicitly and given provisional position
$\widetilde{\vXd}^{w,m}$, the provisional velocities are obtained by
evaluating
\begin{equation}
  \label{eq:sdc-vel-gi}
  \widetilde{\vectord{u}}^{w,m}
  = \linopd{D}_{\widetilde{\gamma}^{w,m}}^{-1}
  \left(\linopd{G}_{\widetilde{\gamma}^{w,m}}
  \vectord{f}(\widetilde{\vXd}^{w,m},\widetilde{\scalard{\sigma}}^{w,m},
  \widetilde{\scalard{\lambda}}^{w,m}) +
  \vectord{u}^{\infty}\right).
\end{equation}
which requires a global inversion of $\linopd{D}_{\widetilde{\gamma}^{w,m}}$.
%
%
The same approach is taken for \li and \lic schemes, except the
provisional velocities are obtained using local inversion only,
all the inter-particle interactions are treated explicitly and added
to the explicit term, i.e., complementary velocity $\widetilde{\comp{\vectord{u}}}_i^{w-1,m}$;
modifying \cref{eq:sdc-vel-gi} for each vesicle, we obtain
\begin{equation}
  \label{eq:sdc-vel-li}
  \widetilde{\vectord{u}}_i^{w,m}
  = \linopd{D}_{\widetilde{\gamma}_i^{w,m}}^{-1}
  \left(\linopd{G}_{\widetilde{\gamma}_i^{w,m}}\vectord{f}(\widetilde{\vXd}_i^{w,m},
  \widetilde{\scalard{\sigma}}^{w,m}_i,\widetilde{\scalard{\lambda}}^{w,m}_i)
  + \widetilde{\comp{\vectord{u}}}_i^{w-1,m}\right),
\end{equation}
where $\widetilde{\comp{\vectord{u}}}_i^{w-1,m}$ is computed using
$\widetilde{\vXd}^{w,m}$ and $\widetilde{\vectord{u}}_j^{w-1,m}$($j \neq i$)
accounting the velocity influence from other vesicles.
We only need to invert the local interaction matrices
$\linopd{D}_{\widetilde{\gamma}_i^{w,m}}$ in this scheme.


\subsubsection{Contact-resolving iteration}
Let $\linopd{A}\vXd^{+}=\vectord{b}$ be the linear system that is
solved at each iteration of a \lic scheme (in case of the
\lic-\abbrev{SDC} scheme, on each of the inner step of the
\abbrev{SDC}).
$\linopd{A}$ is a block diagonal matrix, with blocks
$\linopd{A}_{ii}$ corresponding to the self interactions of the
$\ordinal{i}$ particle.  All inter-particle interactions are treated
explicitly, and thus included in the right-hand side $\vectord{b}$.
We write \cref{eq:vesicle-bie-disc}, or \cref{eq:sdc-error-disc}, in a compact form as
\begin{align}
  \linopd{A}\vXd^{+} &= \vectord{b} + \linopd{G}\vectord{f}_c^{+},
  \label{eq:stepping-colex-vel-short}\\
  0 \le \discrete{V}(\discrete{\gamma};t^{+}) \quad & \bot \quad \discrete{\lambda} \ge 0, \label{eq:stepping-colex-collision-short}
\end{align}
which is a mixed Nonlinear Complementarity Problem (\abbrev{NCP}),
because the \abbrev{STIV} function $\discrete{V}(\discrete{\gamma},t)$
is a nonlinear function of position.
Note that since this is the \lic scheme, $\linopd{G}$ is also a block
diagonal matrix.
To solve this \abbrev{NCP}, we use a first-order linearization of the
$\discrete{V}(\discrete{\gamma};t)$ to obtain an \abbrev{LCP} and
iterate until the \abbrev{NCP} is solved to the desired accuracy:
\begin{align}
  \linopd{A}\vXd^{\star} = \vectord{b} + \linopd{G}\linopd{J}^T\discrete{\lambda}.
  \label{eq:stepping-colex-vel-LCP}\\
  0 \le \discrete{V}(\discrete{\gamma};t^{+k}) + \linopd{J}  \Delta \vXd \quad  \bot \quad
  \discrete{\lambda} \ge 0, \label{eq:stepping-colex-collision-LCP}
\end{align}
where $\Delta \vXd$ is the update to get the new candidate
solution $\vXd^\star$, and $\linopd{J}$ denotes the Jacobian of the volume
$\nabla_X\discrete{V}(\discrete{\gamma},t^{+k})$.
\Cref{alg:colStep} summarizes the steps to solve
\cref{eq:stepping-colex-vel-short,eq:stepping-colex-collision-short}
as a series of linearization steps
\cref{eq:stepping-colex-vel-LCP,eq:stepping-colex-collision-LCP}.
We discuss the details of the \abbrev{LCP} solver separately below.
\IncMargin{1em}
\begin{algorithm}[!b]
\mcaption{alg:colStep}{Contact-free time-stepping}{}
\DontPrintSemicolon
\SetKwInOut{Input}{input}
\SetKwInOut{Output}{output}
\SetKwFunction{LCPSolver}{lcpSolver}
\SetKwFunction{LCPMatrix}{formLCPMatrix}
\SetKwFunction{ColVolume}{getContactVolume}
\SetKwFunction{ColVGrad}{getContactVolumeJacobi}
\SetKwFunction{ApplyLCPMatrix}{applyLCPMatrix}
\Input{$\vXd, \vectord{b}$}
\Output{$\vXd^{+}, \vectord{f}_c^+$}
  $\linopd{A} \leftarrow  \linopd{A}(\vXd)$\;
  $\vectord{b} \leftarrow \vectord{b}(\vXd)$\;
  $\vectord{f}_c^{+} \leftarrow 0$\;
  $k \leftarrow 0$\;
  $\vXd^{\star} \leftarrow \linopd{A}^{-1}\vectord{b}$\;
  $\discrete{V} \leftarrow$ \ColVolume{$\vXd^{\star}$}\;
  \While{$\discrete{V}<0$}{
  $\linopd{J} \leftarrow$ \ColVGrad{$\vXd^{\star}$}\;
  $\lambda \leftarrow$ \LCPSolver{$\discrete{V}$}\;
  $k \leftarrow k+1$\;
  $\vectord{b} \leftarrow \vectord{b} + \linopd{GJ}^T\lambda$\;
  $\vXd^{\star} \leftarrow \linopd{A}^{-1}\vectord{b}$\;
  $\discrete{V} \leftarrow$ \ColVolume{$\vXd^{\star}$}\;
  $\vectord{f}_c^{+} \leftarrow \vectord{f}_c^{+}+\linopd{J}^T\lambda$\;
}
  $\vXd^{+} \leftarrow \vXd^{\star}$\;
\end{algorithm}
\DecMargin{1em}

In lines $1$ to $6$, we solve the unconstrained system
$\linopd{A}\vXd^{\star}=\vectord{b}$
using the solution from previous time step. Then, the
\abbrev{STIV}s are computed to check for collision.
The loop at lines $7{-}14$ is the linearized contact-resolving steps.
Substituting \cref{eq:stepping-colex-vel-LCP} into
\cref{eq:stepping-colex-collision-LCP}, and using the fact that
$\Delta \vXd = \linopd{A}^{-1}\linopd{GJ}^T\lambda$ we cast the
problem in the standard \abbrev{LCP} form
\begin{equation}
  \label{eq:LCP}
  0 \le \discrete{V} + \linopd{B}\discrete{\lambda} \quad \bot
  \quad \discrete{\lambda} \ge 0,
\end{equation}
where $\linopd{B}=\linopd{JA}^{-1}\linopd{GJ}^T$. The \abbrev{LCP}
solver is called on line $9$ to obtain the magnitude of the constraint
force, which is in turn used to update obtain new candidate positions
that may or may not satisfy the constraints.
Line $13$ checks the minimal-separation constraints for the candidate
solution.
In line 11, the collision force is incorporated into the
right-hand-side $\vectord{b}$ for self interaction in the next
\abbrev{LCP} iteration.
In line 14, the contact force is updated, which will be used to form
the right-hand-side $\vectord{b}$ for the global interaction in the
next time step.

The \abbrev{LCP} matrix $\linopd{B}$ is an $M$ by $M$ matrix,
where $M$ is the number of contact volumes, $M = \mathcal{O}(N_v + N_p)$.
Each entry $\linopd{B}_{k,p}$ is induced change in the $\ordinal{k}$
contact volume by the \ordinal{p} contact force.
Matrix $\linopd{B}$ is sparse and typically diagonally dominant, since
most \abbrev{STIV} volumes are spatially separated.

\subsection{Solving the Linear Complementarity Problem}
\label{sec:lcp-solver}
In the contact-resolving iterations we solve an \abbrev{LCP}, \cref{eq:LCP}.
Most common algorithms (e.g., Lemke's algorithm
\cite{lemke1965bimatrix} and splitting based iterative algorithms
\cite{mangasarian1977solution,ahn1981solution})
requires explicitly formed \abbrev{LCP} matrix $\linopd{B}$, which can be
prohibitively expensive when there are many collisions.
We use the minimum-map Newton method, which we  modify to require
matrix-vector evaluation only, as we can perform it without explicitly forming the matrix.

We briefly summarize the minimum-map Newton method.
Let $\vectord{y} = \discrete{V}+\linopd{B} \lambda$. Using the minimum map
reformulation we can convert the \abbrev{LCP} to a root-finding problem
\begin{equation}
  \label{eq:LCP-minmap}
  \vectord{H}(\discrete{\lambda}) \equiv \begin{bmatrix}h(\lambda_1,y_1)
  \\ \cdots \\ h(\lambda_M,y_M) \end{bmatrix} = 0,
\end{equation}
where $h(\lambda_i,y_i) = \min(\lambda_i,y_i)$.
This problem is solved by Newton's method (\cref{alg:lcpSolver}). In the algorithm,
$\linopd{P}_{\mathbb{A}}$ and $\linopd{P}_{\mathbb{F}}$
are selection matrices: $\linopd{P}_{\mathbb{A}}\lambda$ selects the rows of
$\discrete{\lambda}$ whose indices are in set $\mathbb{A}$ and zeros out
all the other rows.
While function $\vectord{H}$ is not smooth, it is Lipschitz and
directionally differentiable, and its so-called B-derivative
$\linopd{P}_{\mathbb{A}}\linopd{B} + \linopd{P}_{\mathbb{F}}$ can be
formed to find the descent direction for Newton's method
\cite{Erleben2013}.
The matrix $\linopd{P}_{\mathbb{A}}\linopd{B} + \linopd{P}_{\mathbb{F}}$ is a
sparse matrix, and we use \gmres to solve this linear system.
Since $\linopd{B}$ is sparse and diagonally dominant, in practice the
linear system is solved in few \gmres iterations and the Newton solver
converges quadratically.

\IncMargin{1em}
\begin{algorithm}[!b]
\mcaption{alg:lcpSolver}{Minimum Map \abbrev{LCP} Solver}{}
\DontPrintSemicolon
\SetKwInOut{Require}{require}
\SetKwInOut{Output}{output}
\SetKw{KwNewton}{Newton's step:\quad}
\SetKw{Tol}{tol}
\SetKwFunction{ProjLineSearch}{projectLineSearch}
\SetKwFunction{ApplyLCPMatrix}{applyLCPMatrix}
  \Require{$\ApplyLCPMatrix()$, $\discrete{V}$ and $\epsilon$}
  \Output{$\lambda$}
  $e \leftarrow \epsilon$\;
  $\lambda \leftarrow 0$\;
\While{$e > \epsilon $}{
  $\vectord{y} \leftarrow \discrete{V} + \text{\ApplyLCPMatrix{$\lambda$}}$\;
  $\mathbb{A} \leftarrow \left\{ i \middle| y_i < \lambda_i \right\}$
  \tcp*[r]{index of active constraints}
  $\mathbb{F} \leftarrow \left\{ i \middle| y_i \ge \lambda_i \right\}$\;
  Iteratively solve $\displaystyle\begin{bmatrix} \linopd{B} &
  \linopd{-I} \\ \linopd{P}_{\mathbb{F}} &
  \linopd{P}_{\mathbb{A}} \end{bmatrix} \begin{bmatrix}\Delta
    \lambda\\\Delta y \end{bmatrix} = \begin{bmatrix}
    0\\ -\linopd{P}_{\mathbb{A}} \vectord{y} - \linopd{P}_{\mathbb{F}}
    \discrete{\lambda}\end{bmatrix}$\tcp*[r]{$\linopd{B}$ applied by applyLCPMatrix}
  $\tau \leftarrow$  \ProjLineSearch{$\Delta \lambda$} \;
  $\lambda \leftarrow \lambda + \tau \Delta \lambda$ \;
  $e \leftarrow \norm{\vectord{H}(\lambda)}$
}
\end{algorithm}
\DecMargin{1em}

\subsection{Algorithm Complexity}
We estimate the complexity of a single time step as a function of
the number of points on each vesicle  $N$, number of vesicles $N_v$.
Let $C_N$ denote the cost of inverting a local linear system for one
particle; then the complexity of inverting linear systems for all
particles is \O{C_N N_v}.
In \cite{Veerapaneni2009,Rahimian2010} it is shown that
for \li scheme $C_N = \O{N\log N}$.
The cost of evaluating the inter-particle interactions
at the $N_v N$ discrete points using \abbrev{FMM} is $\mathcal{O}(N_v N)$.

We assume that for each contact resolving step, the number of contact
volumes is $M$. Assuming that minimum map Newton method takes $K_1$
steps to converge, the cost of solving the \abbrev{LCP} is \O{K_1 C_N
  N_v}, because inverting $\linopd{A}$ is the costliest step in
applying the \abbrev{LCP} matrix.
The total cost of solving the \abbrev{NCP} problem is \O{K_1 K_2 C_N
  N_v}, where $K_2$ are the number of contact resolving iterations.
  In the numerical simulations we observe that the minimum map Newton
  method converges in a few iterations ($K_1\approx 15$) and the
  number of contact resolving iterations is also small and independent
  of the problem size ($K_2\approx 10$).
In \cref{sec:results}, we compare the cost of solving contact
constrained system and the cost of unconstrained system.

\section{Results\label{sec:results}}
In this section, we present results characterizing the accuracy,
robustness, and efficiency of a locally-implicit time stepping scheme
combined with our contact resolution framework in comparison to
schemes with no contact resolution.
\begin{itemize}
  \item First, to demonstrate the robustness of our scheme in
    maintaining the prescribed minimal separation distance with
    different viscosity contrast $\nu$, we consider two vesicles in an
    extensional flow, \cref{ss:extensionalFlow}.
  \item In \cref{ss:shearFlow}, we explore the effect of minimal
    separation $\minsep$ and its effect on collision displacement in
    shear flow.
    We demonstrate that the collision scheme has a minimal effect on
    the shear displacement.
  \item We compare the cost of our scheme with the unconstrained
    system using a simple sedimentation example in \cref{ss:sedFlow}.
    While the per-step cost of the unconstrained system is marginally
    lower, it requires a much spatial and temporal resolution in order
    to maintain a valid contact-free configuration, making the overall
    cost prohibitive.
  \item We report the convergence behavior of different time-stepping
    in \cref{ssc:conv} and show that our scheme achieves second order
    convergence rate with \sdc[2].
  \item We illustrate the efficiency and robustness of our
    algorithm with two examples: 100 sedimenting vesicles in a
    container and a flow with multiple vesicles and rigid particles
    within a constricted tube in \cref{ss:stenosisFlow}.
\end{itemize}


Our experiments support the general observation that when vesicles
become close, the \li scheme does a very poor job in handling of
vesicles' interaction \cite{rahimian2015} and the time stepping
becomes unstable.
The \gi scheme stays stable, but the iterative solver requires more
and more iterations to reach the desired tolerance, which in turn
implies higher computational cost for each time step.
Therefore, the \gi scheme performance degrades to the point of not
being feasible due to the intersection of vesicles and the \li scheme
fails due to intersection or the time-step instability.

\subsection{Extensional flow\label{ss:extensionalFlow}}
To demonstrate the robustness of our \colres framework, we consider
two vesicles placed symmetrically with respect to the $y$ axis in the
extensional flow $\vu = [-x,y]$.
The vesicles have reduced area of $0.9$ and we use a first-order time
stepping with \li, \lic, and \gi schemes for the experiments in this
test.
We run the experiments with different time step size and viscosity
contrast and report the minimal distance between vesicles as well as
the final error in vesicle perimeter, which should be kept constant
due local inextensibility.
Snapshots of the vesicle configuration for two of the time-stepping
schemes are shown in \cref{fig:extenImpVisc256Snap}.

\begin{figure}[!b]
  \centering
  \setlength\figurewidth{1.7in}
  \setlength\figureheight{1.2in}
  \hspace{-28pt}\includepgf{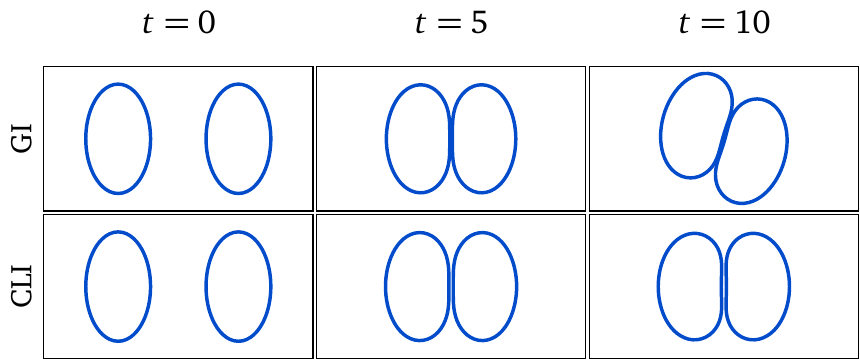}%
  \mcaption{fig:extenImpVisc256Snap}{ Snapshots of two vesicles in
    extensional flow using \gi and \lic schemes}{ As the distance
    between two vesicles decreases the configuration loses symmetry in
    the \gi scheme as shown in the top row.  Nevertheless, as shown in
    the second row, the \li with minimal-separation constraint scheme
    maintains the desired minimum separation distance and two vesicles
    also maintain a symmetric configuration.  (The viscosity contrast
    is 500 in this simulation).  }
\end{figure}

\begin{figure}[!bt]
  \centering
  \setlength\figurewidth{.3\linewidth}
  \setlength\figureheight{.24\linewidth}
  \hspace{20pt}
  \subfloat[The minimal distance of two vesicles over time \label{sfg:exten-dis-time}]{\hspace{-40pt}\includepgf{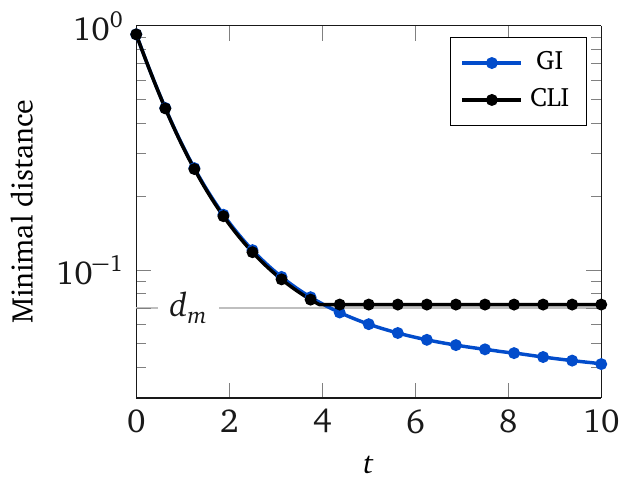}}
  \hspace{40pt}
  \subfloat[The minimum distance between two vesicles versus the viscosity contrast\label{sfg:exten-dis-visc}]{\hspace{-20pt}\includepgf{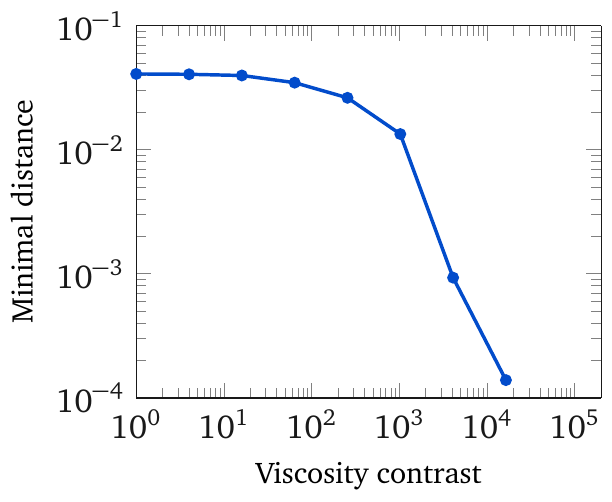}}
  \mcaption{fig:exten-dis}{Distance between two vesicles in extensional flow}{
    \subref{sfg:exten-dis-time} The minimum distance between two
    vesicles over time for both \lic and \gi schemes.
    The \lic scheme easily maintains the prescribed minimal separation
    of $h$.
    \subref{sfg:exten-dis-visc} The final distance (at $T=10$) between
    two vesicles as viscosity contrast is increased in \gi scheme.
  }
\end{figure}

In \cref{sfg:exten-dis-time}, we plot the minimal distance between two
vesicles over time. The vesicles continue to get closer in the \gi
scheme. However, the \lic scheme maintains the desired minimum
separation distance between two vesicles.
In \cref{sfg:exten-dis-visc}, we show the minimum distance between the
vesicles over the course of simulation (with time horizon $T=10$)
versus the viscosity contrast.
As expected, we observe that the minimum distance between two vesicles
decreases as the viscosity contrast is increased.
Consequently, for higher viscosity contrast with both \gi and \li
schemes, either the configuration loses its symmetry or the two
vesicles intersect.
With minimal-separation constraint, we any desired minimum separation
distance between vesicles is maintained, and the simulation is more
robust and accurate as shown in \cref{fig:extenImpVisc256Snap} and
\cref{tbl:extenViscError}.

\begin{table}[!t]
  \newcolumntype{C}{>{\centering\let\newline\\\arraybackslash\hspace{0pt}$}c<{$}}
  \newcolumntype{R}{>{\centering\let\newline\\\arraybackslash\hspace{0pt}$}r<{$}}
  \newcolumntype{L}{>{\centering\let\newline\\\arraybackslash\hspace{0pt}$}l<{$}}
  \newcolumntype{A}[1]{>{\centering\let\newline\\\arraybackslash\hspace{0pt}}D{.}{.}{#1}}
  \setlength{\tabcolsep}{4pt}
  \newcommand\cc[2]{#1}
  \newcommand\ce[3]{\cc{#1e{-#2}}{#3}}
  \centering
  \small
  \begin{tabular}{LLLLL}
    \toprule
    \nu & \dt    & \text{\lic}       & \text{\li}        & \text{\gi}        \\ \cmidrule(lr){1-1}\cmidrule(lr){2-2}\cmidrule(lr){3-5}
    1   & 0.4    & \ce{1.17}{01}{6}  & \ce{1.17}{01}{6}  & \ce{7.66}{02}{12} \\
    1   & 0.2    & \ce{9.49}{04}{24} & \ce{9.49}{04}{24} & \ce{1.03}{03}{18} \\
    1   & 0.1    & \ce{4.49}{04}{24} & \ce{4.49}{04}{24} & \ce{4.69}{04}{24} \\
    1   & 0.05   & \ce{2.23}{04}{24} & \ce{2.23}{04}{24} & \ce{2.29}{04}{24} \\
    1   & 0.025  & \ce{1.12}{04}{24} & \ce{1.12}{04}{24} & \ce{1.14}{04}{24} \\
    1   & 0.0125 & \ce{5.65}{05}{30} & \ce{5.65}{05}{30} & \ce{5.68}{05}{30} \\ \cmidrule(lr){1-1}\cmidrule(lr){2-2}\cmidrule(lr){3-5}
    1e2  & 0.4    & \ce{9.42}{03}{18} & -                 & \ce{2.54}{04}{24} \\
    1e2  & 0.2    & \ce{1.33}{04}{24} & \ce{1.33}{04}{24} & \ce{1.22}{04}{24} \\
    1e2  & 0.1    & \ce{6.38}{05}{30} & \ce{6.38}{05}{30} & \ce{5.96}{05}{30} \\
    1e2  & 0.05   & \ce{3.05}{05}{30} & \ce{3.05}{05}{30} & \ce{2.95}{05}{30} \\
    1e2  & 0.025  & \ce{1.49}{05}{30} & \ce{1.49}{05}{30} & \ce{1.47}{05}{30} \\
    1e2  & 0.0125 & \ce{7.39}{06}{36} & \ce{7.39}{05}{30} & \ce{7.32}{06}{36} \\
    \bottomrule
  \end{tabular}
  \hspace{12pt} 
  \begin{tabular}{LLLLL}
    \toprule
   \nu  & \dt    & \text{\lic}               & \text{\li}        & \text{\gi}        \\ \cmidrule(lr){1-1}\cmidrule(lr){2-2}\cmidrule(lr){3-5}
    1e3 & 0.4    & \ce{8.88}{03}{18}^{T=8.8} & -                 & \ce{1.66}{05}{30} \\
    1e3 & 0.2    & \ce{2.08}{02}{12}         & -                 & \ce{7.99}{06}{36} \\
    1e3 & 0.1    & \ce{5.85}{06}{36}         & -                 & \ce{3.93}{06}{36} \\
    1e3 & 0.05   & \ce{2.42}{06}{36}         & \ce{1.93}{06}{36} & \ce{1.95}{06}{36} \\
    1e3 & 0.025  & \ce{1.12}{06}{36}         & \ce{9.78}{07}{42} & \ce{7.33}{07}{42} \\
    1e3 & 0.0125 & \ce{5.89}{07}{42}         & \ce{4.87}{07}{42} & \ce{2.01}{07}{42} \\ \cmidrule(lr){1-1}\cmidrule(lr){2-2}\cmidrule(lr){3-5}
    1e4 & 0.4    & \ce{1.45}{03}{18}^{T=3.6} & -                 & -                 \\
    1e4 & 0.2    & \ce{7.23}{04}{24}^{T=5.4} & -                 & -                 \\
    1e4 & 0.1    & \ce{7.75}{04}{24}^{T=6.8} & -                 & -                 \\
    1e4 & 0.05   & \ce{2.41}{06}{36}         & -                 & -                 \\
    1e4 & 0.025  & \ce{1.17}{06}{36}         & -                 & -                 \\
    1e4 & 0.0125 & \ce{5.58}{07}{42}         & -                 & -                 \\
    \bottomrule
  \end{tabular}

  \mcaption{tbl:extenViscError}{Error in the length of vesicles in
    extensional flow}{The error in the final length of two vesicles in
    extensional flow with respect to viscosity contrast, timestep
    size, and for different schemes.
    The experiment's setup is described in \cref{ss:extensionalFlow}
    and snapshots of which are shown in \cref{fig:extenImpVisc256Snap}.
    The cases with a ``$-$'' indicate that either vesicles have
    intersected or the \gmres failed to converge due to
    ill-conditioning of the system; the latter happens in the \gi
    scheme and high viscosity contrast.
    Cases with superscript indicate that the flow loses its symmetry
    at that time.
  }
\end{table}

In \cref{tbl:extenViscError}, we report the final error in vesicle
perimeter for different schemes with respect to viscosity contrast and
timestep size.
With minimal-separation constraint, we achieve similar or smaller
error in length compared to \li or \gi methods (when these methods
produce a valid result).
Moreover, one can use relatively large time step in all flow
parameters\emdash/most notably when vesicles have high viscosity
contrast.
In contrast, the \li scheme requires very small, often impractical,
time steps to prevent instability or intersection.

\subsection{Shear flow\label{ss:shearFlow}}
We consider vesicles and rigid bodies in an unbounded shear flow and
explore the effects of minimal separation on shear diffusivity.
In the first simulation, we consider two vesicles of reduced area
$0.98$ (to minimize the effect of vesicles' relative orientation on
the dynamics) placed in a shear flow with shear rate $\chi = 2$.
Let $\delta_t = |y_{t}^{1} - y_{t}^{2}|$ denote the vertical offset
between the centroids of vesicles at time $t$.
Initially, two vesicles are placed with a relative vertical offset
$\delta_0$ as show in \cref{fig:shear-snap}.
\begin{figure}[!bt]
  \centering
  \subfloat[$t=5$\label{sfg:shear-snap1}]{\framebox(96,50){\includegraphics[angle=0,width=.2\linewidth]{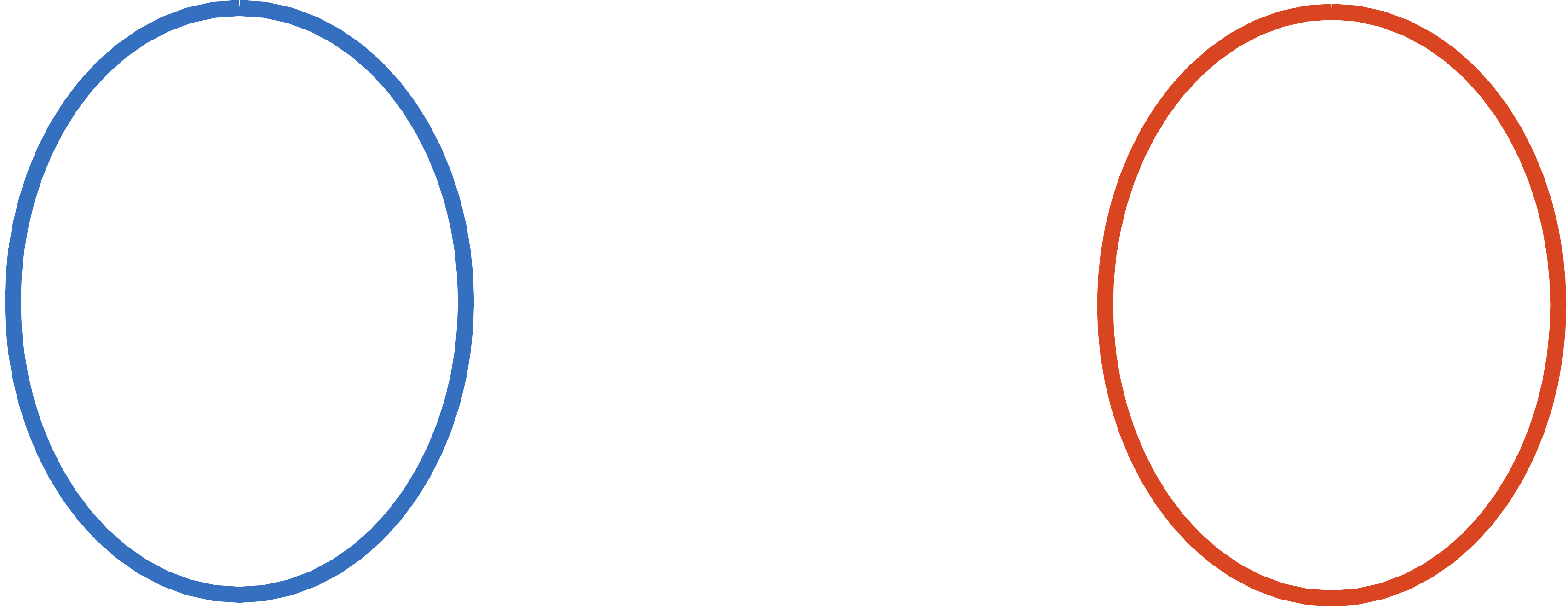}}}
  \hspace{4pt}
  \subfloat[$t=7.5$\label{sfg:shear-snap2}]{\framebox(91,50){\includegraphics[angle=0,width=.19\linewidth]{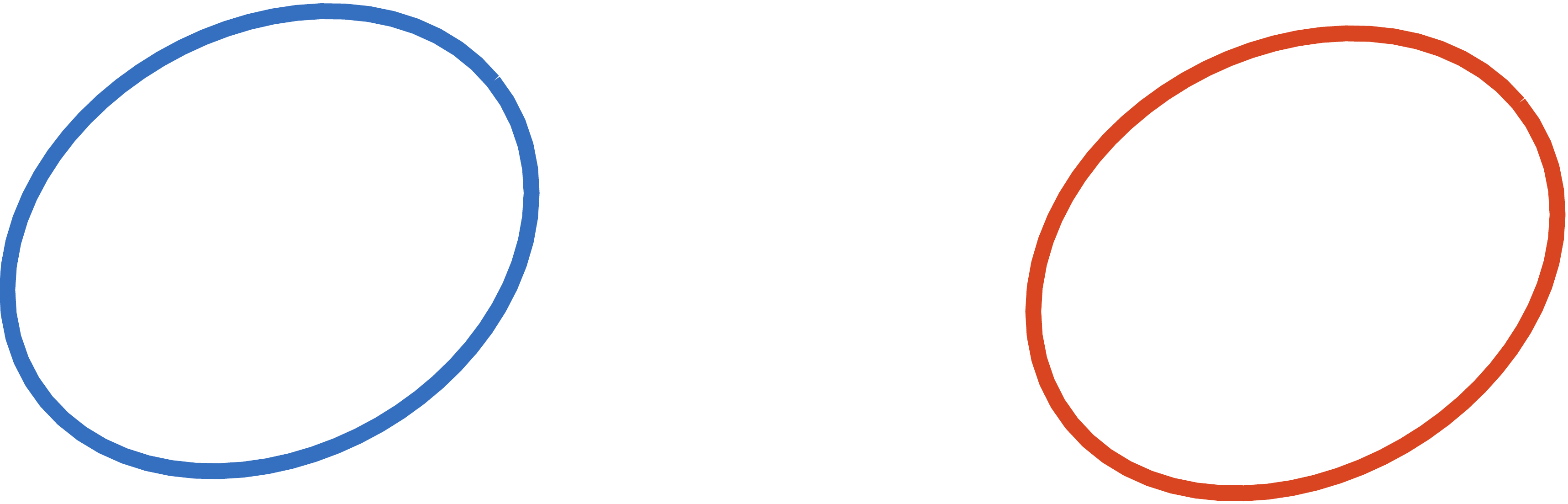}}}
  \hspace{4pt}
  \subfloat[$t=10$\label{sfg:shear-snap3}]{\framebox(60,50){\includegraphics[angle=0,width=.09\linewidth]{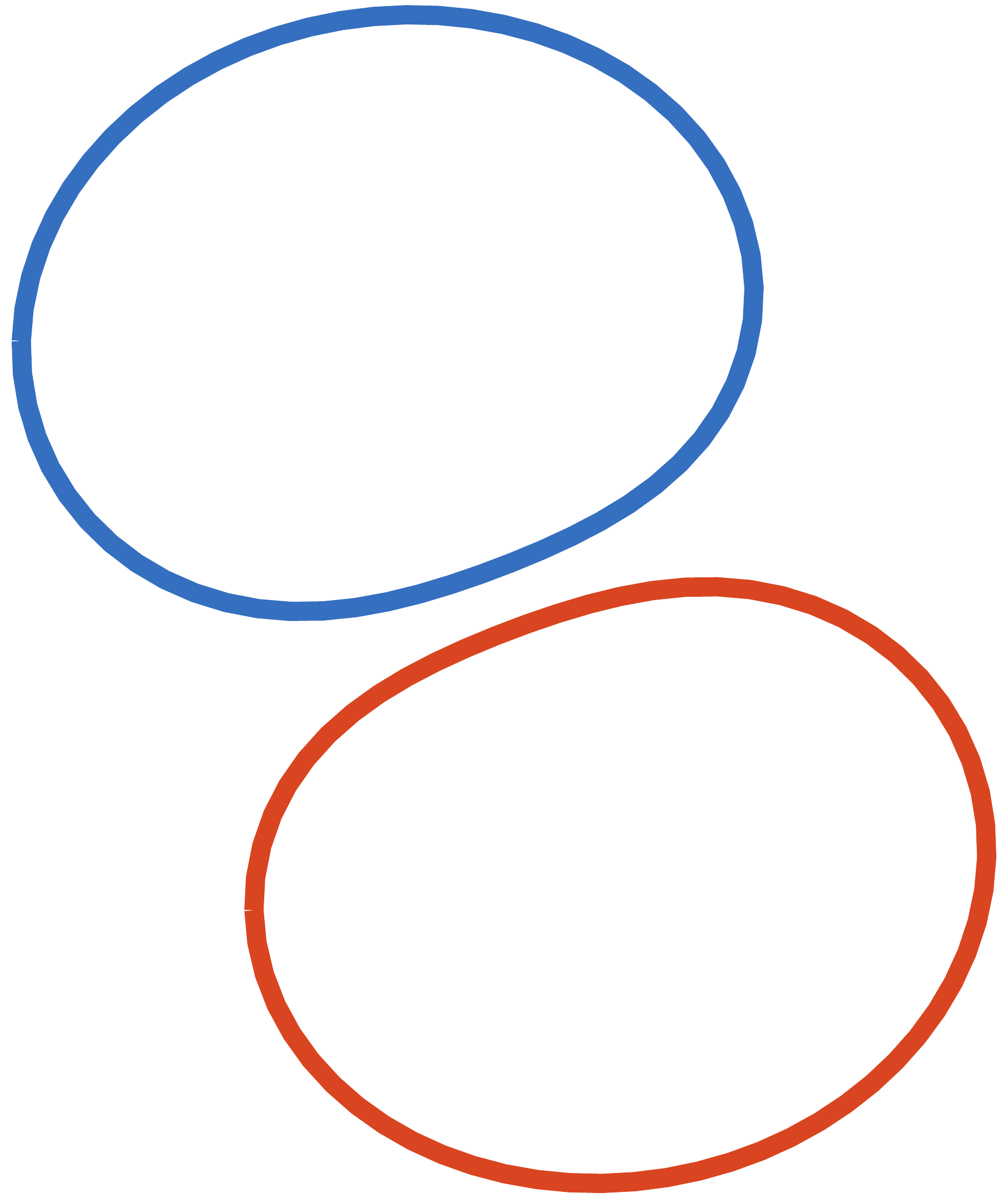}}}
  \hspace{4pt}
  \subfloat[$t=12$\label{sfg:shear-snap4}]{\framebox(77,50){\includegraphics[angle=0,width=.16\linewidth]{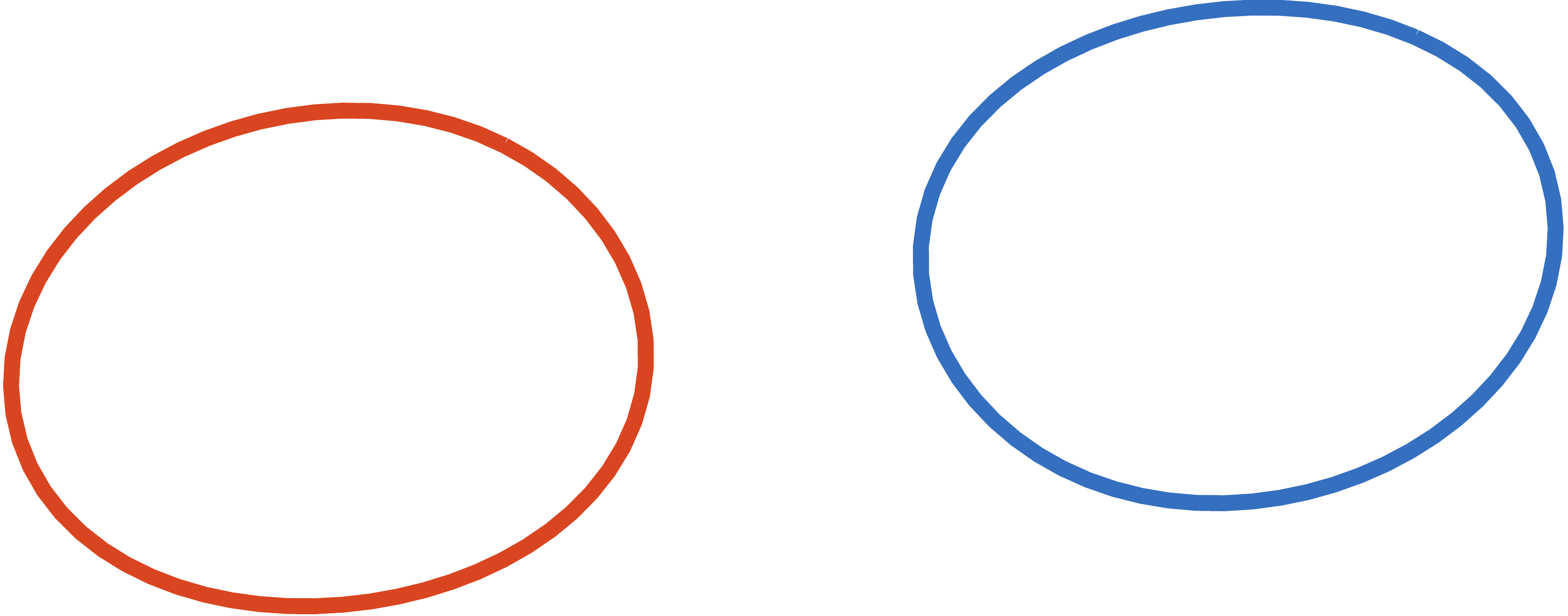}}}
  \mcaption{fig:shear-snap}{Shear flow experiment setup}{The snapshots
    of two vesicles in shear flow.  Initially, one vesicle is placed
    at $[-8, \delta_0]$ and the second vesicle is placed at $[0, 0]$.}
\end{figure}

\begin{figure}[!bt]
  \setlength\figurewidth{.24\linewidth}
  \setlength\figureheight{.2\linewidth}
  \hspace{.5em}
  \subfloat[Centroid trajectory (vesicle)\label{sfg:shear-cm}]{\hspace{-30pt}\includepgf{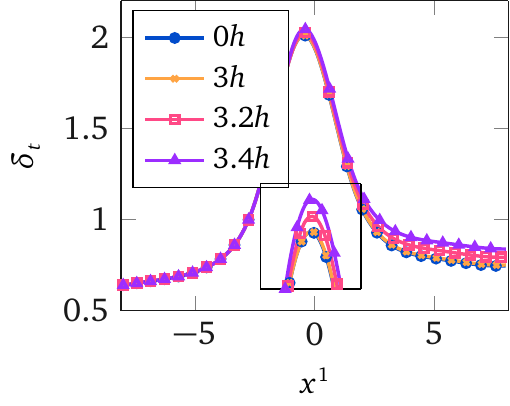}}
  \hspace{1em}
  \subfloat[Centroid trajectory (rigid)\label{sfg:shear-cmrgd}]{\hspace{-30pt}\includepgf{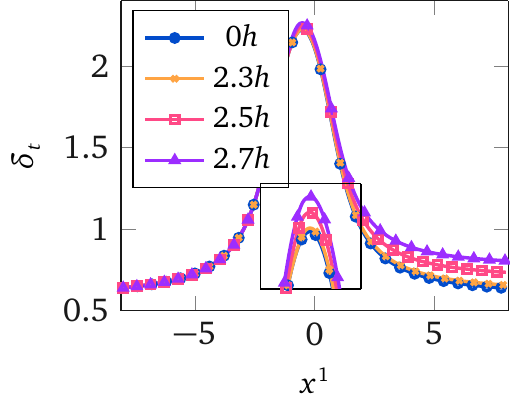}}
  \hspace{1em}
  \subfloat[Final centroid offset\label{sfg:shear-sd}]{\hspace{-20pt}\includepgf{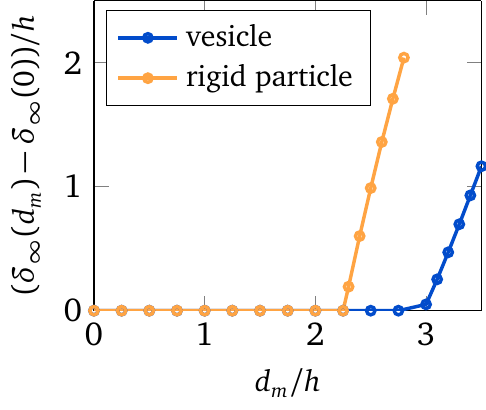}}
\mcaption{fig:cenMass}{The offset $\delta_t(\minsep)$ between the
  centroids of two vesicles in shear flow}{The initial offset is
  $\delta_0 = 0.64$ and $N=64$ discretization points are used,
  implying $h=0.0994$, where $h$ is the distance between two
  discretization points along vesicle boundary.  \subref{sfg:shear-cm}
  and \subref{sfg:shear-cmrgd} show $\delta_t(x^1)$ for the vesicles
  and circular rigid particles for different minimal separation
  distance $\minsep$.
  \subref{sfg:shear-sd} shows the excess displacement induced by
  minimal separation. When collision constraint is activate (larger
  $\minsep$), the particles are in effect hydrodynamically larger and
  the excess displacement grows linearly with $\minsep$.
}
\end{figure}

In \cref{fig:cenMass}, we report $\delta_{t}$ and its value upon
termination of the simulation when $x^1>8$, denoted by
$\delta_{\infty}(\minsep)$.
In \cref{sfg:shear-cm,sfg:shear-cmrgd}, we plot $\delta_{t}$ with
respect to the $x^1$ for different minimum separation distances.
Based on a high-resolution simulation~(with $N=128$ and $8\times$
smaller time step), the minimal distance between two vesicles without
contact constraint is about $2.9h$ for vesicles and $2.2h$ for rigid
particles.
As the minimum separation parameter $\minsep$ is decreased below this
threshold, the simulations with minimal-separation constraint converge
to the reference simulation without minimal-separation constraint.
In \cref{sfg:shear-sd}, we plot the excess terminal displacement due
to contact constraint, $(\delta_{\infty}(d) - \delta_{\infty}(0))/h$,
as a function of the minimum separation distance. When collision
constraint is activate, the particles are in effect hydrodynamically
larger and the excess displacement grows linearly with $\minsep$.

\subsection{Sedimentation\label{ss:sedFlow}}
To compare the performance of schemes with and without contact
constraints, we first consider a small problem with three vesicles
sedimenting in a container.
We compare  three first-order time stepping schemes: locally
implicit~(\li), locally implicit with collision handling~(\lic), and
globally implicit~(\gi).
\begin{figure}[!bt]
  \centering
  \setlength\figureheight{.3\linewidth}
  \subfloat[\lic $t=0$\label{sfg:sed3nv-snap1}]{\includegraphics[height=\figureheight]{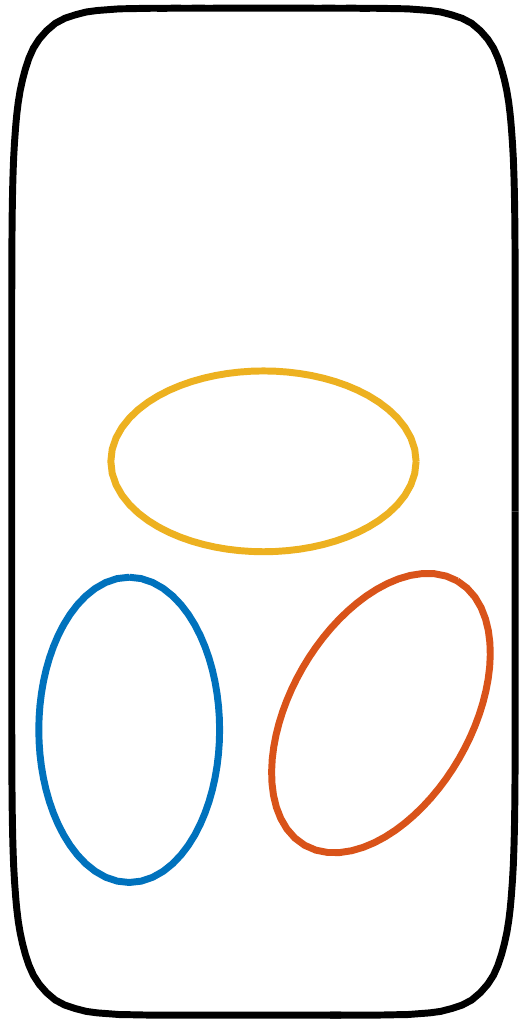}}
  \subfloat[\lic $t=5$\label{sfg:sed3nv-snap2}]{\includegraphics[height=\figureheight]{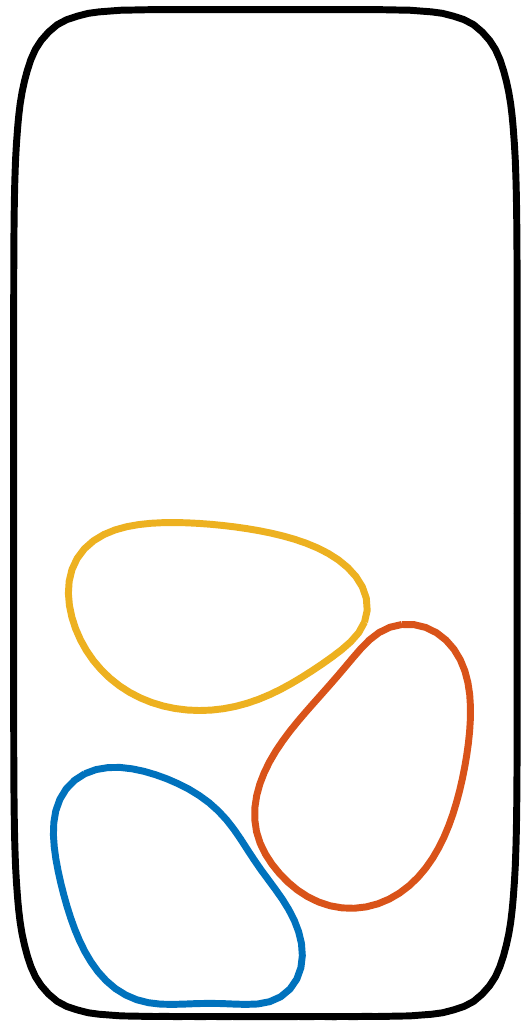}}
  \subfloat[\mbox{\lic $t=20$\hspace{-2pt}}\label{sfg:sed3nv-snap3}]{\includegraphics[height=\figureheight]{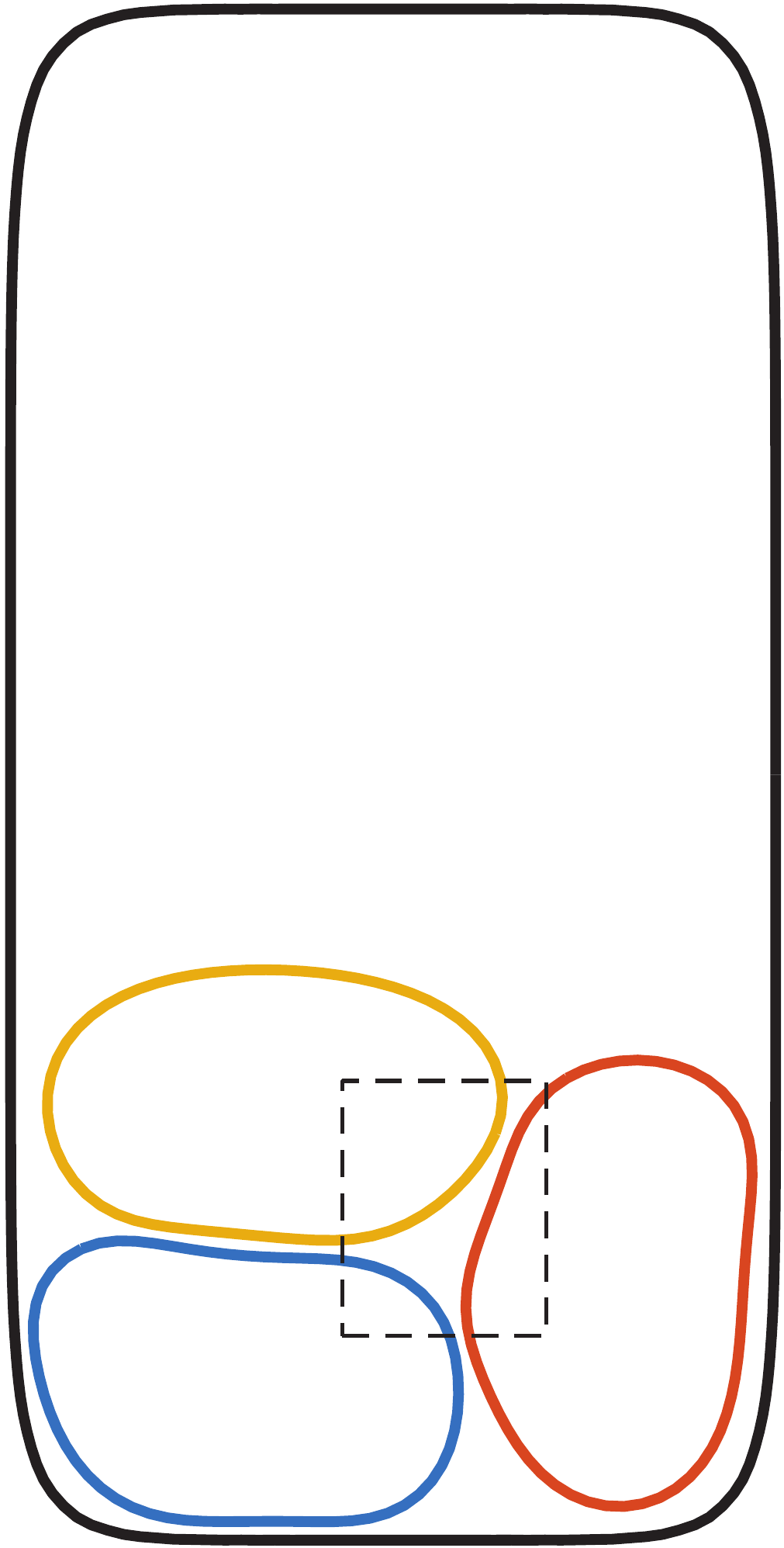}}
  \subfloat[\lic $t=20$\label{sfg:sed3nv-snap3z}]{
    \tikzsetnextfilename{external-sed3nv-snap3z}%
    \tikz{
      \node (A1) {\includegraphics[height=.7\figureheight]{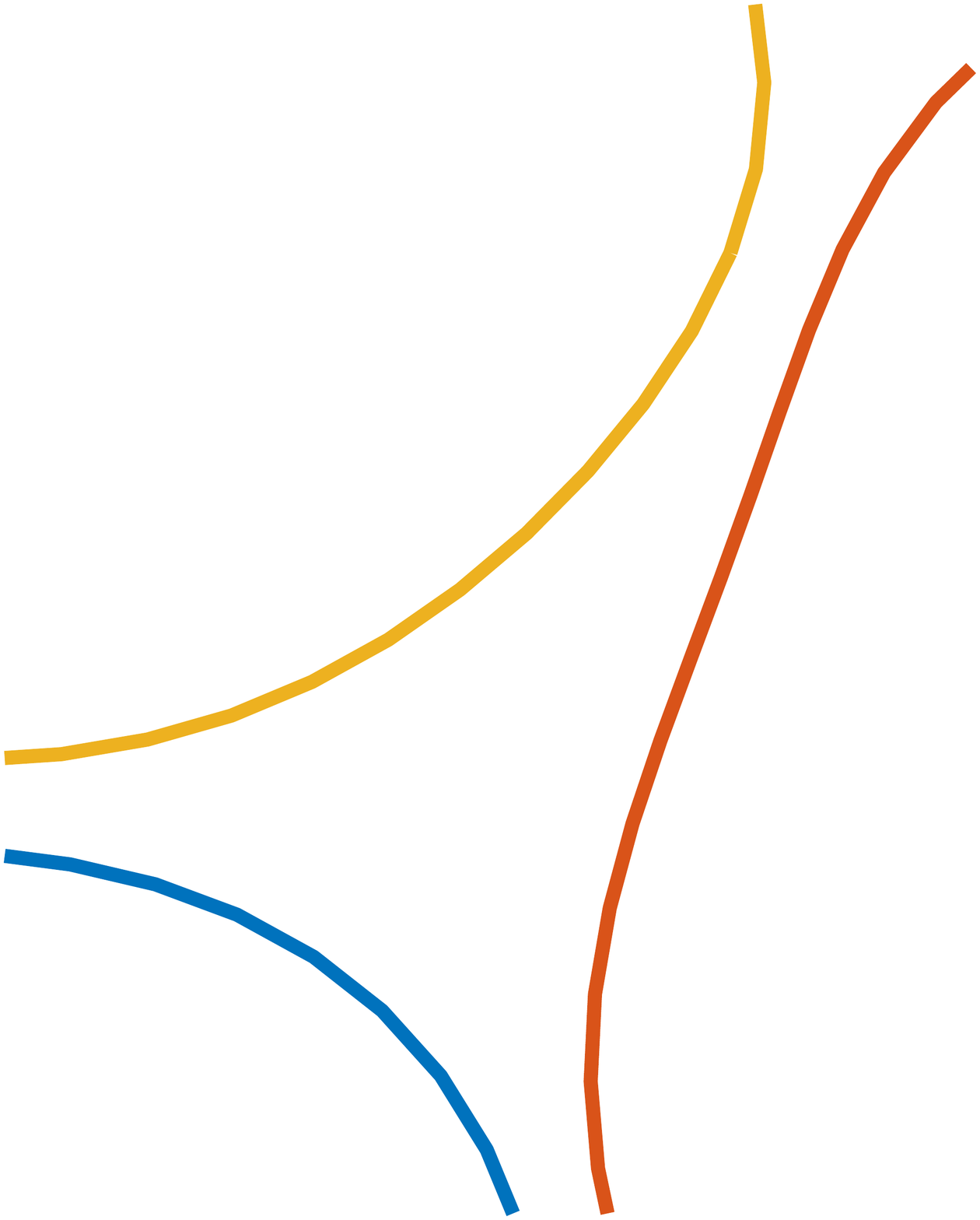}};
      \node [draw, dashed, fit=(A1)]{};
    }
  }
  \subfloat[\gi $t=10$\label{sfg:sed3nv-impcol}]{\includegraphics[height=\figureheight]{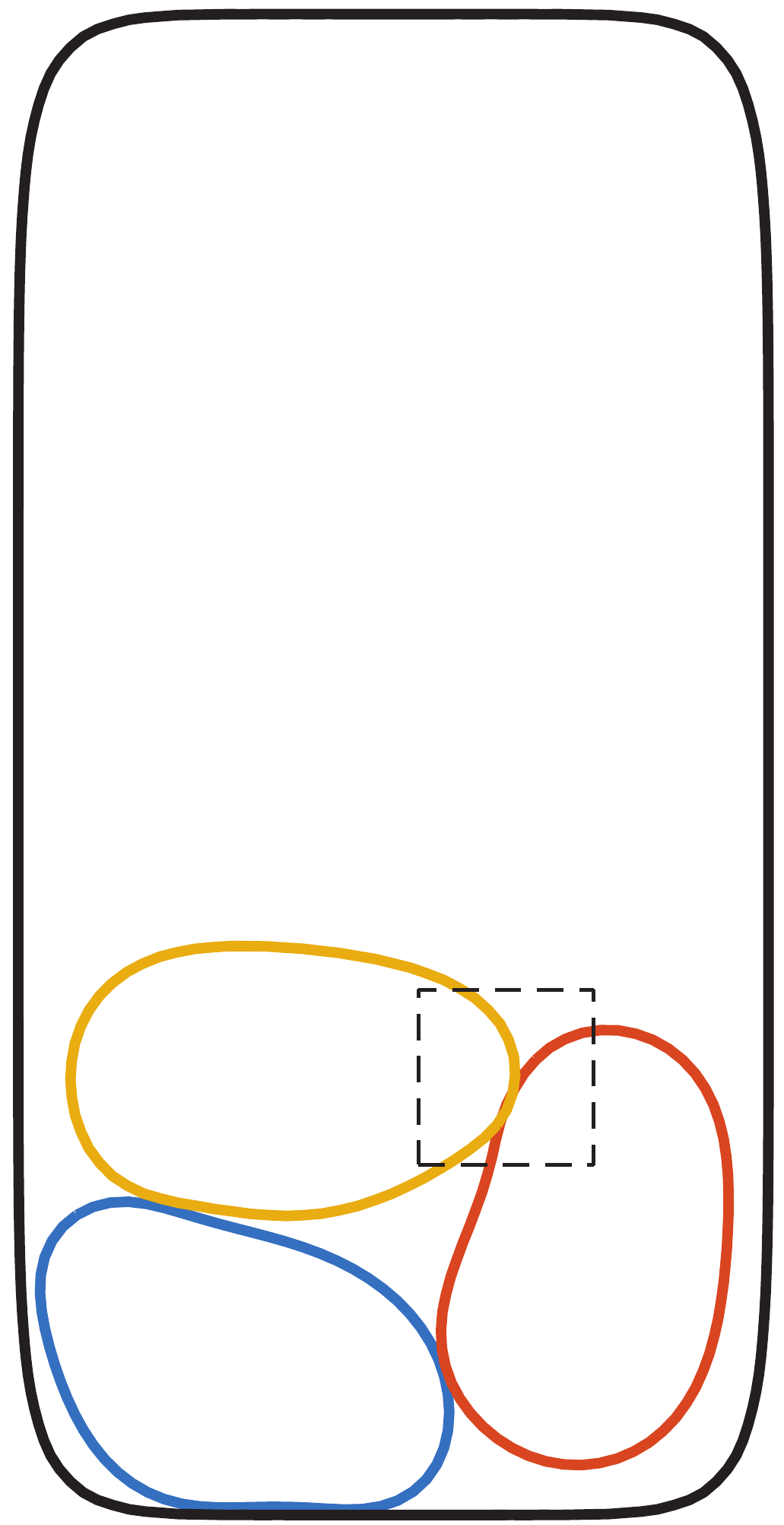}}
  \subfloat[\gi $t=10$\label{sfg:sed3nv-impcolz}]{
    \tikzsetnextfilename{external-sed3nv-impcolz}%
    \tikz{
      \node (B1) {\includegraphics[height=.7\figureheight]{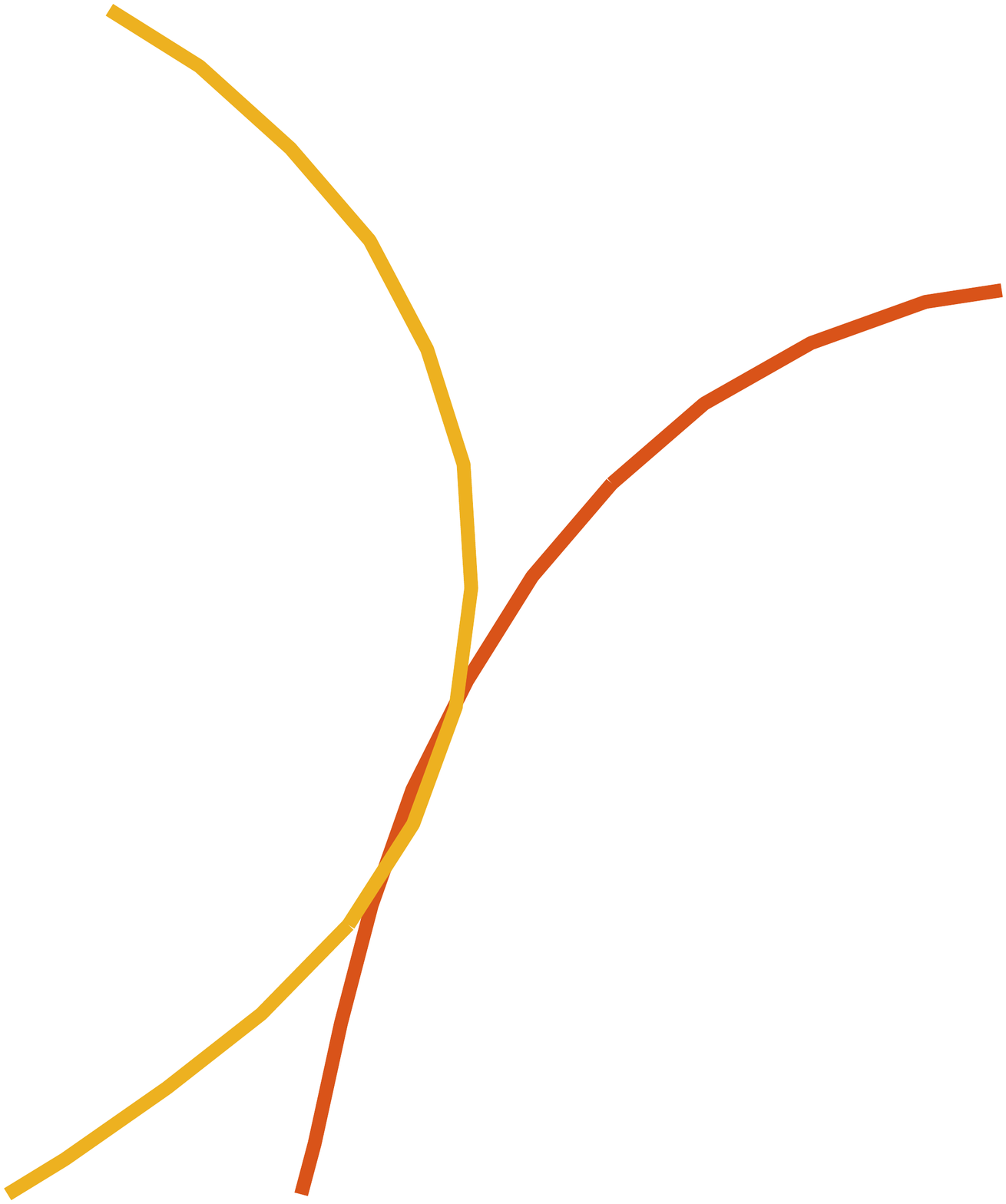}};
      \node [draw, dashed, fit=(B1)]{};
    }
  }
  \mcaption{fig:sed3nv-snap}{Sedimenting vesicles}{Snapshots of three
    sedimenting vesicles in a container using three time-stepping
    schemes.
    The vesicles have reduced area of $0.9$, their viscosity contrast
    is $100$, and are discretized with $64$ points.  The boundary is
    discretized with $256$ points, the simulation time horizon $T=26$,
    and the time step is $0.01$.
    \subref{sfg:sed3nv-snap1}--\subref{sfg:sed3nv-snap3} The
    dynamics of three vesicles sedimenting with \lic scheme.
    \subref{sfg:sed3nv-snap3z} The contact region of
    \subref{sfg:sed3nv-snap3}.  \subref{sfg:sed3nv-impcol} The
    instance when the vesicles intersect using \gi scheme.
    \subref{sfg:sed3nv-impcolz} The contact region of
    \subref{sfg:sed3nv-impcol}.}
\end{figure}

Snapshots of these simulations are shown in \cref{fig:sed3nv-snap}.
For the \li scheme, the error grows rapidly when the vesicles
intersect and a $64\times$ smaller time step is required for resolving
the contact and for stability.
Similarly, for the \gi scheme the vesicle intersect as shown in
\cref{sfg:sed3nv-impcol,sfg:sed3nv-impcolz} and a $4\times$ smaller
time step is needed for the contact to be resolved.
The \lic scheme, maintains the desired minimal separation between
vesicles.

Although the current code is not optimized for computational
efficiency, it is instructive to consider the relative amount of time
spent for a single time step in each scheme.
Each time step in \li scheme takes about $1$ second, the time goes up
to $1.5$ seconds for the \lic scheme.
In contrast, a single time step with the \gi scheme takes, on average,
$65$ seconds because the solver needs up to 240 \gmres iterations to
converge when vesicles are very close.
This excessive cost renders this scheme impractical for large
problems.

To demonstrate the capabilities of the \lic scheme and to gain a
qualitative understanding of the scaling of the cost as the number of
intersections increases, we consider the sedimentation of $100$
vesicles.
As the sediment progresses the number of contact regions grows to
about $70$ per time step.
For this simulation, we use a lower viscosity contrast of $10$ and the
enclosing boundary is discretized with $512$ points and the total time
is $T = 30$.
Snapshots of the simulation with \lic scheme are shown in
\cref{fig:sed-snap}.

We observe that with first-order time stepping, both \li and \lic
schemes are unstable.
Therefore, we run this simulation with second order \sdc using \li and
\lic schemes.
The \gi scheme is prohibitively expensive and infeasible in this case,
due to ill-conditioning and large number of \gmres iterations per
time-step.

The \li scheme requires at least $12000$ thousand time steps to
maintain the non-intersection constraint and each time step takes
about $700$ seconds to complete.
On the other hand, the \lic scheme only need $1500$ time steps to
complete the simulation (taking the stable time step) and each time
step takes about $830$ seconds.
We repeated this experiment with $16$ discretization points on each
vesicle using \lic scheme and $1500$ time steps are also sufficient
for this case (with final length error about $3.83\e{-3}$), each time
step takes about $170$ seconds and the simulation takes about $70$
hours to complete. The number of contact resolving iterations in
\cref{alg:colStep} is about $10$.

In summary, \li scheme requires $4\times$ more points on each vesicle
and $8\times$ smaller time-step size to keep vesicles in a valid
configuration compared to \lic scheme.

\subsection{Convergence analysis\label{ssc:conv}}
\begin{figure}[!bt]
  \centering
  \setlength\figurewidth{3.5in}
  \setlength\figureheight{1.5in}
  \hspace*{-28pt}\includepgf{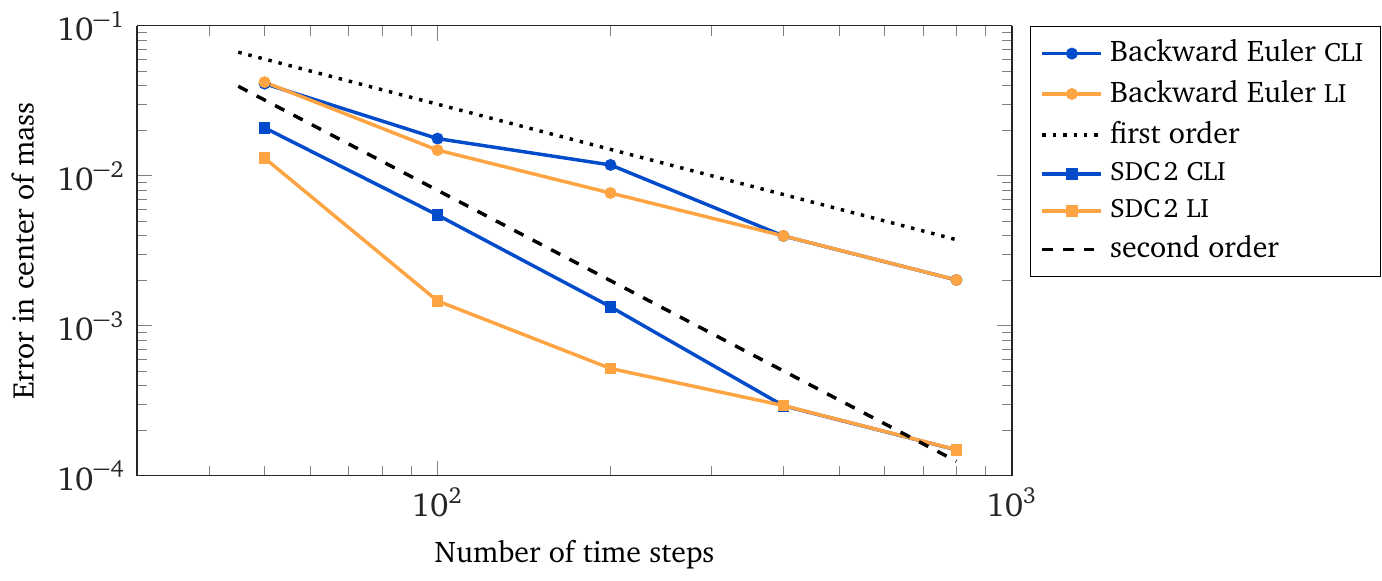}%
  \mcaption{fig:conv}{Convergence rate}{We compare the final center of
    mass of three sedimenting vesicles in a container for
    Backward Euler and \sdc[2] (second order Spectral Deferred Correction
    scheme).
    The time horizon is set to $T=2$. We choose the spatial resolution
    proportional to the time-step.
    As a reference, we use the results for the \gi scheme with time
    step $\dt = 1.25\e{-3}$ and $N=256$ discretization points on each
    vesicle and $512$ discretization points on the boundary.
  }
  \vspace{-4pt}
\end{figure}
To investigate the accuracy of the time-stepping schemes, we consider
the sedimentation of three vesicles with reduced area $0.9$ in a
container as shown in \cref{fig:sed3nv-snap}.
We test \li and \lic schemes for a range of time steps and spatial
resolutions and report the error in the location of the center of mass
of the vesicles at the end of simulation.
The spatial resolution $h$ is chosen proportional to time-step size
and for the \lic scheme the minimal separation $\minsep$ is
proportional to $h$.
As a reference, we use a fine-resolution simulation with \gi scheme.
The error as a function of time step size is shown in \cref{fig:conv}.

\subsection{Stenosis\label{ss:stenosisFlow}}
\begin{figure}[!bt]
  \centering
  \setlength\figurewidth{.7\linewidth}
  \subfloat[$t=0$ \label{sfg:choke-snap1}]{\includegraphics[angle=0,width=\figurewidth]{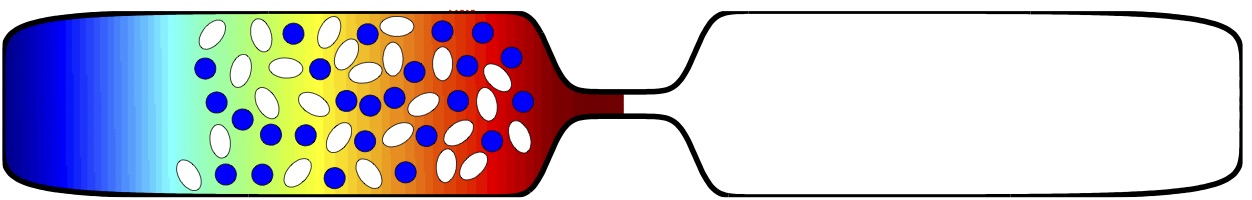}}\\
  \subfloat[$t=4.5$ \label{sfg:choke-snap2}]{\includegraphics[angle=0,width=\figurewidth]{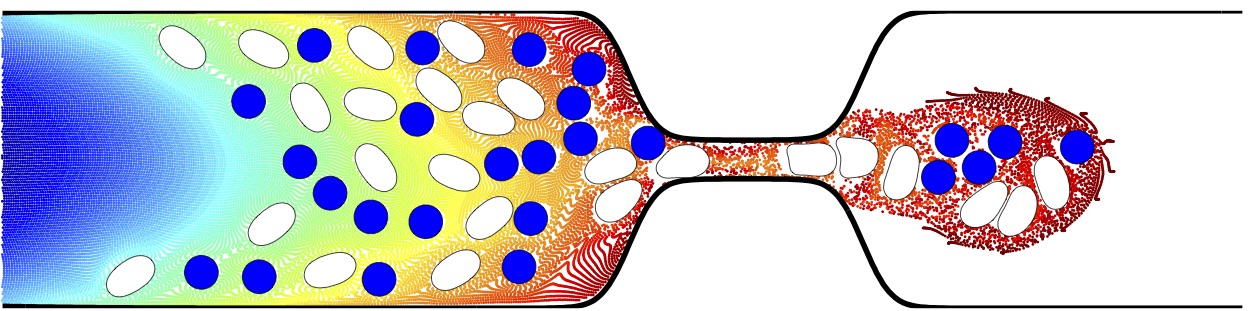}}\\
  \subfloat[$t=6.1$ \label{sfg:choke-snap3}]{\includegraphics[angle=0,width=\figurewidth]{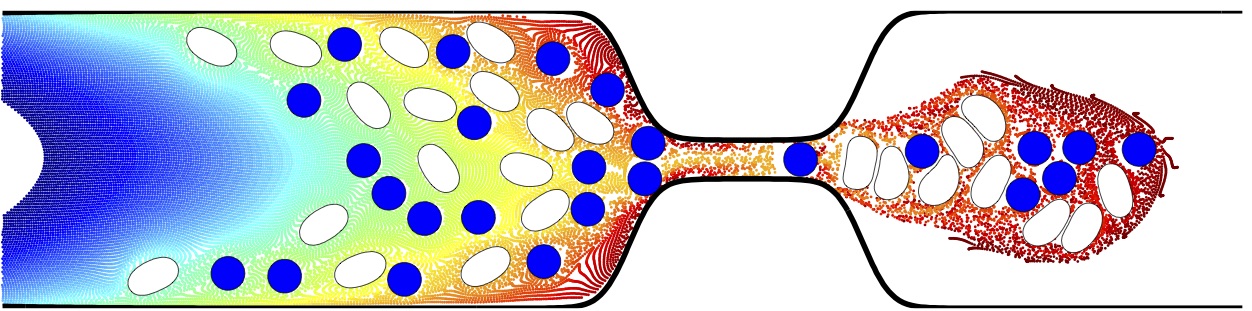}}\\
  \subfloat[$t=11.3$ \label{sfg:choke-snap4}]{\includegraphics[angle=0,width=\figurewidth]{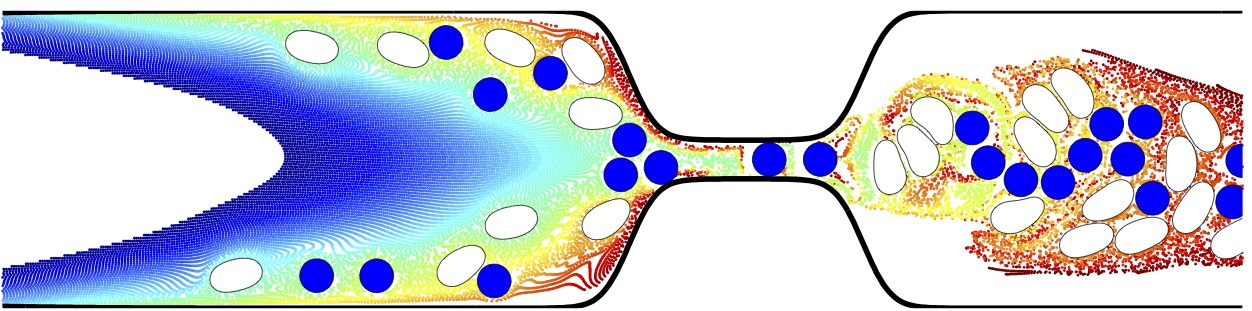}}
  \mcaption{fig:choke-snap}{Stenosis with $25$ vesicles and $25$ rigid
    particles}{Snapshots of vesicles and rigid particles passing a
    constricted tube, the fluid flows from left to right (using proper
    Dirichlet boundary condition on the wall).
    The colored tracers are for visualization purposes and do not
    induce any flow.
    The simulation time horizon is set to $T=20$, each vesicle or
    rigid particle is discretized with $32$ points, and the wall is
    discretized with $1024$ points.
    \subref{sfg:choke-snap1} The initial configuration.
    \subref{sfg:choke-snap2}--\subref{sfg:choke-snap4}
    The interaction and collision between vesicles, rigid particles, and
    the domain boundary at different instances.
    Without minimal-separation constraint, vesicles and rigid
    particles easily contact at the entrance of the constricted area.
  }
\end{figure}
\begin{figure}[!bt]
  \centering
  \setlength\figurewidth{3.4in}
  \setlength\figureheight{1.8in}
  \hspace*{-20pt}\includepgf{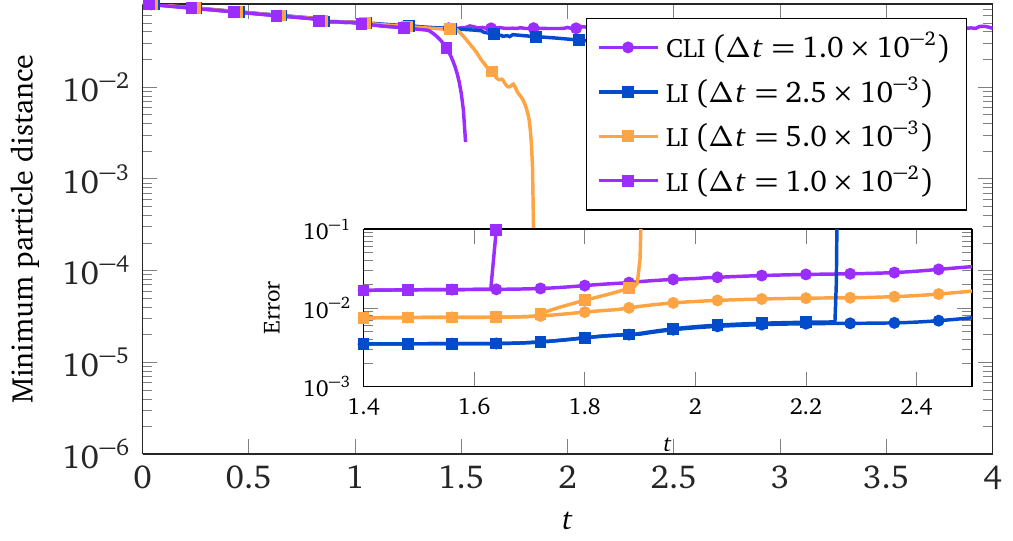}%
  \mcaption{fig:choke-error}{error for stenosis example}{The minimal
    distance and error for different schemes and time step sizes.  For
    the $\li$ scheme, independently of the time step size, vesicles or
    particles intersect with the domain boundary, resulting in
    divergence.  The \lic scheme is stable and maintains the desired
    minimal distance.  }
\end{figure}

As a stress test for particle-boundary and
vesicle-boundary interaction, we study the flow with $25$ vesicles of
reduced area $0.9$, mixed with $25$ circular rigid particles in a
constricted tube (\cref{fig:choke-snap}).
It is well known that rigid bodies can become arbitrarily close in the
Stokes flow, e.g. \cite{loewenberg1997}, and without proper collision
handling, the required \andor{temporal}{spatial} resolution would be
very high.

In this example, the vesicles and rigid particles are placed at the
right hand side of the constricted tube.
We use backward Euler method and search for the stable time step for
schemes \li and \lic. Similarly to the sedimentation example, we do
not consider the $\gi$ scheme due to its excessive computational cost.

The largest stable time step for the \lic scheme is $\dt = 0.01$.
For the \li scheme, we tested the cases with up to $32\times$ smaller
time step size, but we were unable to avoid contact and intersection
between vesicles and rigid particles.

\Cref{fig:choke-error} shows the error and minimal distance between
vesicles, particles, and boundary for \lic and \li schemes with
different time step sizes.
Without the minimal-separation constraint, the solution diverges when
the particles cross the domain boundary.

To validate our estimates for the errors due to using a
piecewise-linear approximation in the minimal separation constraint
instead of high-order shapes, we plot the minimum distance at each
step for the piecewise-linear approximation and upsampled shapes in
\cref{fig:chokediscf}.
The target minimal separation distance is set to $\minsep=h$.
We observe that the actual minimal distance for smooth curve is
smaller than the minimal distance for piecewise-linear curve, while
the difference between two distances is small compared to the target
minimum separation distance.
\begin{figure}[!bt]
  \centering
  \setlength\figurewidth{2.5in}
  \setlength\figureheight{1.5in}
  \hspace*{-28pt}\includepgf{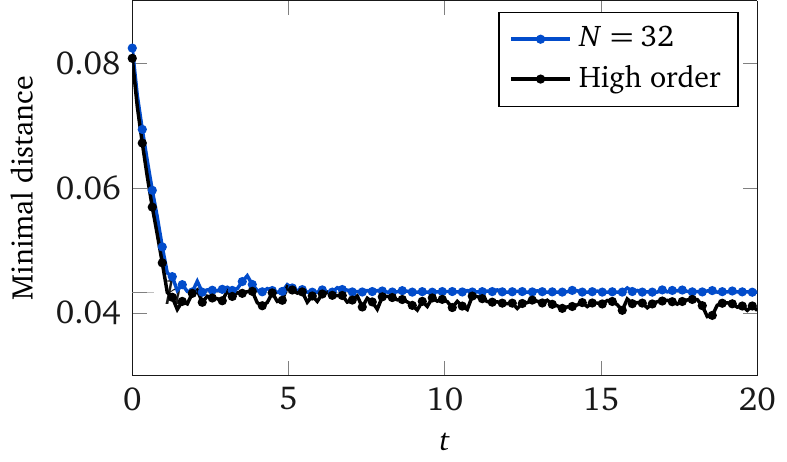}
  \mcaption{fig:chokediscf}{Effective minimal distance}{The minimal distance
    between 32-point piecewise-linear approximations (enforced by the
    constraint) and the true minimal distance between high-order
    shapes.
    We use an upsampled linear approximation with $N=128$ as the
    surrogate for the high-order curves.
  }
\end{figure}

Note that the due to higher shear rate in the constricted area, the
stable time step is dictated by the dynamics in that area.  For such
flows, we expect a combination of adaptive time stepping
\cite{Quaife2014a} and the scheme outlined in this paper provide the
highest speedup.

\section{Conclusion\label{sec:conclusion}}
In this paper we introduced a new scheme for efficient simulation of
dense suspensions of deformable and rigid particles immersed in
Stokesian fluid.
We demonstrated through numerical experiments that this scheme is
orders of magnitude faster than the alternatives and can achieve high
order temporal accuracy.
We are working on extending this approach to three dimensions and
using adaptive time stepping.

We extend our thanks to George Biros, David Harmon, Ehssan Nazockdast,
Bryan Quaife, Michael Shelley, and Etienne Vouga for stimulating
conversations about various aspects of this work.

\appendix

\printbibliography 

\end{document}